\documentclass{article}
\usepackage{iclr2026_conference,times}
\usepackage{hyperref}
\newcommand{\titlemeta}{FastLSQ: Solving PDEs in One Shot via Fourier Features with Exact Analytical Derivatives}

\newcommand{\authormeta}{Antonin Sulc\orcidlink{0000-0001-7767-778X}\\
Lawrence Berkeley National Laboratory, Berkeley, U.S.A.\\
\texttt{asulc@lbl.gov}}

\newcommand{\abstractmeta}{We present FastLSQ, a framework for PDE solving and inverse problems built on trigonometric random Fourier features with exact analytical derivatives.
Trigonometric features admit closed-form derivatives of any order in $\mathcal{O}(1)$, enabling graph-free operator assembly without autodiff.
Linear PDEs: one least-squares call; nonlinear: Newton--Raphson reusing analytical assembly.
On 17 PDEs (1--6D), FastLSQ achieves $10^{-7}$ in 0.07\,s (linear) and $10^{-8}$--$10^{-9}$ in $<$9\,s (nonlinear), orders of magnitude faster and more accurate than iterative PINNs.
Analytical higher-order derivatives yield a differentiable digital twin; we demonstrate inverse problems (heat-source, coil recovery) and PDE discovery.
Code: github.com/sulcantonin/FastLSQ and \texttt{pip install fastlsq}.}

\usepackage{amssymb}
\usepackage{amsmath}
\usepackage{amsthm}
\usepackage{mathtools}
\usepackage{graphicx}
\usepackage{booktabs}
\usepackage{multirow}
\usepackage{xcolor}
\usepackage{algorithm}
\usepackage{algpseudocode}
\usepackage{hyperref}
\usepackage{orcidlink}
\usepackage{tikz}
\usepackage{xspace}
\usepackage{enumitem}
\usepackage{wrapfig}
\usepackage{subcaption}
\usepackage{listings}
\iclrfinalcopy
\usetikzlibrary{arrows.meta, positioning, shapes.geometric, calc, backgrounds, fit}

\newtheorem{theorem}{Theorem}[section]

\newtheorem{proposition}[theorem]{Proposition}
\newtheorem{corollary}[theorem]{Corollary}

\newcommand{\RR}{\mathbb{R}}
\newcommand{\NN}{\mathbb{N}}

\newcommand{\calL}{\mathcal{L}}

\newcommand{\calN}{\mathcal{N}}
\newcommand{\calU}{\mathcal{U}}
\newcommand{\calO}{\mathcal{O}}
\newcommand{\calB}{\mathcal{B}}
\newcommand{\bx}{\mathbf{x}}
\newcommand{\bW}{\mathbf{W}}
\newcommand{\bA}{\mathbf{A}}

\newcommand{\bb}{\mathbf{b}}
\newcommand{\bbeta}{\boldsymbol{\beta}}
\newcommand{\norm}[1]{\left\| #1 \right\|}

\newcommand{\fastlsq}{\textsc{Fast-LSQ}\xspace}

\selectfont
\title{\titlemeta{}}

\author{\authormeta{}}

\begin{document}

\maketitle

\begin{abstract}
\abstractmeta{}
\end{abstract}
\section{Introduction}
\label{sec:introduction}

Solving partial differential equations (PDEs) is a cornerstone of computational physics, with applications from fluid dynamics and electromagnetism to quantum mechanics and climate modeling.
Classical numerical methods, finite differences, finite elements, and spectral methods, remain the workhorses of scientific computing, but face challenges in high dimensions and can require substantial problem-specific implementation effort.

The emergence of physics-informed neural networks (PINNS)~\citep{raissi2019physics} offered a mesh-free alternative that parametrizes the PDE solution as a neural network and minimizes a physics-based loss via stochastic gradient descent.
While conceptually elegant, PINNs require minutes to hours of iterative training \citep{hao2023pinnacle}, suffer from spectral bias \citep{tancik2020fourier}, causality violations \citep{wang2022respecting}, and sensitivity to loss balancing.

Random feature methods (RFMs) offer a middle ground.
By approximating the solution as a linear combination of randomly sampled basis functions with frozen parameters, RFMs reduce the problem to fitting a linear output layer.
For linear PDEs, this yields a linear system in the output coefficients for any choice of basis, as noted by PIELM \citep{dwivedi2020pielm} (using $\tanh$) and RF-PDE \citep{liao2024solving}.
However, RF-PDE still requires between 600 and 2000 epochs of iterative optimization.
One-shot methods like PIELM use $\tanh$ activations, which lack a closed-form cyclic derivative structure; PIELM therefore relies on manual, problem-specific symbolic calculus to assemble $\calL[\tanh]$ for each PDE operator.

\vspace{-0.75em}\paragraph{Key observation.}
While any frozen-feature model yields a linear system for linear PDEs, the practical bottleneck is assembling the PDE operator matrix $A_{ij} = \calL[\phi_j](\bx_i)$.
For sinusoidal features $\phi_j(\bx) = \sin(\bW_j \cdot \bx + b_j)$, we show that $\calL[\phi_j]$ admits an exact closed-form expression for any linear differential operator of any order, requiring only $\calO(1)$ operations per entry and no automatic differentiation.
This is a consequence of the cyclic derivative structure of sinusoids: the $n$-th derivative cycles through $\{\sin, \cos, -\sin, -\cos\}$ with a monomial weight prefactor.
Alternative bases such as $\tanh$ lack any comparable closed-form pattern (Proposition~\ref{prop:tanh}), so methods like PIELM require manual, problem-specific derivation of $\calL[\tanh]$ for each new PDE operator; there is no universal formula.
\fastlsq's operator-agnostic advantage is that one formula~\eqref{eq:cyclic} applies to any $\calL$.

\begin{figure}[t]
\centering
\resizebox{\linewidth}{!}{
\begin{tikzpicture}[
    box/.style={rectangle, draw, rounded corners, minimum width=2.4cm, minimum height=0.55cm, align=center, font=\footnotesize},
    arrow/.style={->, >=stealth, thick},
    label/.style={font=\scriptsize\itshape, text=gray},
]
\node[font=\footnotesize\bfseries, anchor=west] at (-6.5, 1.7) {PINNs / RF-PDE:};
\node[box, fill=red!10] (nn) at (-3.5, 1.05) {Neural Net / RFM};
\node[box, fill=red!10] (ad) at (0.0, 1.05) {Autodiff $\nabla_\theta$};
\node[box, fill=red!10] (loss) at (3.2, 1.05) {$\|\calL[u]-f\|^2$};
\node[box, fill=red!10] (opt) at (6.5, 1.05) {SGD / NLS};
\draw[arrow] (nn) -- (ad);
\draw[arrow] (ad) -- (loss);
\draw[arrow] (loss) -- (opt);
\draw[arrow, dashed, red!60] (opt.north) to[out=30,in=150, looseness=0.18] node[above=-2pt, label] {600--2000 iters} (nn.north);

\node[font=\footnotesize\bfseries, anchor=west] at (-6.5, 0.1) {PIELM (tanh):};
\node[box, fill=orange!10] (tanhbasis) at (-3.5, -0.45) {$\psi_j = \tanh(\bW_j\!\cdot\!\bx\!+\!b_j)$};
\node[box, fill=orange!10] (tanhad) at (1.5, -0.45) {Handwritten $\calL[\psi_j]$};
\node[box, fill=orange!15, draw=orange!50!black] (tanhsolve) at (6.5, -0.45) {$\bbeta = \texttt{lstsq}(\bA, \bb)$};
\draw[arrow] (tanhbasis) -- (tanhad);
\draw[arrow, orange!50!black] (tanhad) -- (tanhsolve);

\node[font=\footnotesize\bfseries, anchor=west] at (-6.5, -1.2) {\textcolor{blue!70!black}{Fast-LSQ (Ours):}};
\node[box, fill=blue!10] (basis) at (-3.5, -1.85) {$\phi_j = \sin(\bW_j\!\cdot\!\bx\!+\!b_j)$};
\node[box, fill=green!15, draw=green!50!black] (solve) at (6.5, -1.85) {$\bbeta = \texttt{lstsq}(\bA, \bb)$};
\draw[arrow, green!50!black, line width=1.0pt] (basis) -- (solve) node[midway, above, label] {closed-form $\calL[\phi_j]$ via \S\ref{sec:derivatives}};
\end{tikzpicture}}
\caption{Architecture comparison. PINNs and RF-PDE require iterative optimization (top row). PIELM and \fastlsq both solve a linear system in one shot, but PIELM requires handwritten, problem-specific calculus to assemble $\calL[\tanh]$ for each PDE operator (middle row), whereas \fastlsq uses a single operator-agnostic formula for any $\calL$ (bottom row).}
\label{fig:architecture}
\vskip -0.1in
\end{figure}

Figure~\ref{fig:architecture} illustrates the key architectural differences.
Both PIELM and \fastlsq are one-shot solvers that reduce to a single least-squares call, so neither involves iterative training.
The distinction lies in how the operator matrix $\bA$ is assembled: PIELM requires handwritten, problem-specific calculus for each PDE operator because $\tanh$ lacks a closed-form pattern, while \fastlsq uses a single operator-agnostic formula~\eqref{eq:cyclic} (\S\ref{sec:derivatives}).
The practical consequences are: (1)~\textbf{operator-agnostic} assembly---one formula for any $\calL$, no manual derivation when changing PDEs; (2)~\textbf{superior spectral accuracy}, with $10\times$ to $1000\times$ lower errors at equal feature counts (\S\ref{sec:experiments}); (3)~\textbf{high-fidelity gradient accuracy}---because derivatives of a sinusoidal basis are just phase-shifted sinusoids scaled by $\calO(\sigma)$, gradient errors remain tightly bounded (within one order of magnitude of value errors; \S\ref{sec:gradient_accuracy}); (4)~\textbf{learnable bandwidth}---$\sigma$ can be made differentiable for inverse problems (Appendix~\ref{app:learnable}); and (5)~a rich set of \textbf{trigonometric symmetries} (eigenfunction property, orthogonality, product-to-sum) underpins both linear and nonlinear PDE solving (Appendix~\ref{app:symmetries}).

We make three contributions: First, we observe that the elementary cyclic derivative structure of sinusoids~\eqref{eq:cyclic} enables graph-free, closed-form assembly of the PDE operator matrix for arbitrary linear differential operators, and we build on this to construct \fastlsq: a one-shot solver that assembles the operator matrix analytically and obtains the solution via a single least-squares call (\S\ref{sec:solver_mode}). Second, we extend \fastlsq to nonlinear PDEs via Newton--Raphson iteration, where each linearized step reuses the analytical assembly, achieving $10^{-8}$ to $10^{-9}$ accuracy in under 9\,s (\S\ref{sec:newton}). Third, on 17 PDEs (5 linear, 5 nonlinear solver-mode, 7 nonlinear regression) in up to 6 dimensions, \fastlsq achieves $10^{-7}$ in under 0.1\,s on linear problems and $10^{-8}$ on nonlinear problems, outperforming all baselines by large margins (\S\ref{sec:experiments}). Furthermore, we run a comparison with RBF Kansa, scikit-fem P2 FEM, and scipy.solve\_bvp on all 10 PDEs, demonstrating that \fastlsq is the only method applicable to every problem, including high-dimensional ($d{\ge}5$) and hyperbolic space-time equations for which conventional solvers have no applicable path.

\section{Method}
\label{sec:method}

\vspace{-0.75em}\subsection{Random Fourier Feature Approximation}
\label{sec:rff}

We approximate the PDE solution $u: \Omega \subset \RR^d \to \RR$ by
\begin{equation}
\label{eq:approximation}
u_N(\bx) = \frac{1}{\sqrt{N}}\sum_{j=1}^{N} \beta_j \, \phi_j(\bx), \qquad
\phi_j(\bx) = \sin(\bW_j^\top \bx + b_j),
\end{equation}
where $\bW_j \sim \calN(\mathbf{0}, \sigma^2 \mathbf{I}_d)$ and $b_j \sim \calU(0, 2\pi)$ are frozen at initialization and only the linear coefficients $\bbeta = (\beta_1, \ldots, \beta_N)^\top$ are trainable.
The $1/\sqrt{N}$ prefactor ensures that the empirical kernel $\frac{1}{N} \sum_{j=1}^{N} \phi_j(\bx)\phi_j(\bx')$ converges to the Gaussian RBF kernel $k(\bx, \bx') = \exp(-\frac{\sigma^2}{2}\norm{\bx - \bx'}^2)$ as $N \to \infty$ \citep{rahimi2007random}, keeping the coefficients $\beta_j$ at $\calO(1)$ magnitude and preventing the ill-conditioning that arises when unnormalized features force $\beta_j \sim 10^6$--$10^8$.
The bandwidth $\sigma$ controls frequency content and the associated kernel length scale.
To capture multiple scales, we use a multi-block architecture with $B$ blocks at potentially different bandwidths $\sigma_b$, concatenating all features into a single vector solved simultaneously.

\vspace{-0.75em}\subsection{Exact Analytical Derivatives of Sinusoidal Features}
\label{sec:derivatives}

The central structural property exploited by \fastlsq is the following.
The multivariate form follows immediately from the chain rule and the cyclic derivatives of $\sin$. For any multi-index $\alpha = (\alpha_1, \ldots, \alpha_d) \in \NN_0^d$ with $|\alpha| = \sum_k \alpha_k$:
\begin{equation}
\label{eq:cyclic}
D^\alpha \phi_j(\bx) = \left(\prod_{k=1}^{d} W_{jk}^{\alpha_k}\right) \cdot \Phi_{|\alpha| \bmod 4}(\bW_j^\top \bx + b_j),
\end{equation}
where $\Phi_0 = \sin$, $\Phi_1 = \cos$, $\Phi_2 = -\sin$, $\Phi_3 = -\cos$.
The identity itself is elementary; what matters is its consequence for PDE solving: for any linear differential operator $\calL$ of any order, every entry of the operator matrix $A_{ij} = \calL[\phi_j](\bx_i)$ reduces to a single trigonometric evaluation multiplied by a monomial in the weights $\bW_j$, computable without automatic differentiation or computational graph construction.
This analytical formulation also enables a differentiable trick: the bandwidth $\sigma$ (or an anisotropic covariance) can be optimised via gradient descent, since gradients flow through the exact inner solve (Appendix~\ref{app:learnable}).

\begin{corollary}[Common operators]
\label{cor:operators}
The Laplacian satisfies $\Delta \phi_j = -\norm{\bW_j}^2 \sin(\bW_j^\top \bx + b_j)$; the biharmonic operator gives $\Delta^2 \phi_j = \norm{\bW_j}^4 \sin(\bW_j^\top \bx + b_j)$; and the advection operator yields $\mathbf{v} \cdot \nabla \phi_j = (\mathbf{v} \cdot \bW_j) \cos(\bW_j^\top \bx + b_j)$.
Each of these is a single function evaluation with a monomial prefactor.
\end{corollary}

\begin{proposition}[Contrast with $\tanh$]
\label{prop:tanh}
The derivative $\frac{d^n}{dz^n}\sin(z) = \Phi_{n\bmod 4}(z)$ is a single trigonometric evaluation for all $n$. No analogous closed-form pattern exists for $\tanh$: while $\frac{d^n}{dz^n}\tanh(z)$ can be evaluated via polynomial recurrences in $\tanh(z)$, the result is a polynomial of degree $n{+}1$ whose coefficients must be tracked, and no fixed-pattern formula reduces it to $\calO(1)$ as in the sinusoidal case.
\end{proposition}

The practical consequence is that assembling the operator matrix entry $\calL[\phi_j](\bx_i)$ for sinusoidal features requires no automatic differentiation and no computational graph:~\eqref{eq:cyclic} reduces it to a single trigonometric evaluation multiplied by a product of weights, computable in $\calO(1)$ regardless of PDE order.
For $\tanh$ features, one must either invoke automatic differentiation, which requires constructing and traversing a computational graph at every collocation point, or implement problem-specific derivative recurrences for each PDE operator.
The sinusoidal closed form avoids both: it is operator-agnostic, graph-free, and numerically stable by construction.

\vspace{-0.75em}\subsection{Solver Mode: Direct Linear Solve}
\label{sec:solver_mode}

For a linear PDE $\calL[u](\bx) = f(\bx)$ in $\Omega$ with boundary conditions $\calB[u](\bx) = g(\bx)$ on $\partial\Omega$, substituting~\eqref{eq:approximation} gives the augmented linear system
\begin{equation}
\label{eq:linear_system}
\underbrace{\begin{pmatrix} \bA^{\text{pde}} \\ \lambda \bA^{\text{bc}} \end{pmatrix}}_{\bA \in \RR^{M \times N}} \bbeta = \underbrace{\begin{pmatrix} \mathbf{f} \\ \lambda \mathbf{g} \end{pmatrix}}_{\bb \in \RR^{M}},
\end{equation}
where $A^{\text{pde}}_{ij} = \calL[\phi_j](\bx_i^{\text{int}})$ is computed in closed form via~\eqref{eq:cyclic}, $A^{\text{bc}}_{ij} = \calB[\phi_j](\bx_i^{\text{bc}})$, and $\lambda > 0$ is a boundary penalty weight.
The system is solved in a single call: $\bbeta^* = \arg\min_{\bbeta} \norm{\bA\bbeta - \bb}_2^2 = \bA^\dagger \bb$, via QR or SVD factorization.

This linearity in $\bbeta$ holds for any frozen-feature basis, including PIELM's $\tanh$.
The advantages of sinusoids are: (1)~operator-agnostic assembly---one formula~\eqref{eq:cyclic} for any $\calL$, versus PIELM's manual derivation per operator; (2)~superior spectral accuracy for smooth and oscillatory solutions (\S\ref{sec:experiments}); and (3)~high-fidelity gradients, because derivatives are phase-shifted sinusoids scaled by $\calO(\sigma)$, so gradient errors stay within one order of magnitude of value errors (\S\ref{sec:gradient_accuracy}).
Algorithm~\ref{alg:solver} summarizes the complete procedure.

\begin{algorithm}[t]
\caption{\fastlsq for Linear PDEs}
\label{alg:solver}
\begin{algorithmic}[1]
\Require PDE operator $\calL$, BC operator $\calB$, source $f$, BC data $g$, domain $\Omega$
\Require Features $N$, bandwidth $\sigma$, collocation counts $M_1, M_2$, penalty $\lambda$
\State Sample frozen weights: $\bW_j \sim \calN(0, \sigma^2 I_d)$, $b_j \sim \calU(0, 2\pi)$
\State Sample collocation points: $\{\bx_i^{\text{int}}\} \subset \Omega$, $\{\bx_i^{\text{bc}}\} \subset \partial\Omega$
\State Build PDE matrix: $A^{\text{pde}}_{ij} \gets \calL[\phi_j](\bx_i^{\text{int}})$ via \S\ref{sec:derivatives} \Comment{Closed-form, no autodiff}
\State Build BC matrix: $A^{\text{bc}}_{ij} \gets \calB[\phi_j](\bx_i^{\text{bc}})$
\State Assemble $\bA, \bb$ per~\eqref{eq:linear_system}; solve $\bbeta^* \gets \texttt{lstsq}(\bA, \bb)$
\State \Return $u_N(\bx) = \frac{1}{\sqrt{N}}\sum_j \beta_j^* \sin(\bW_j^\top \bx + b_j)$
\end{algorithmic}
\end{algorithm}

\vspace{-0.75em}\subsection{Newton--Raphson Extension for Nonlinear PDEs}
\label{sec:newton}

For a nonlinear PDE $\calL[u] + \calN[u] = f$ where $\calN$ is a nonlinear operator, the residual $\calL[u_N] + \calN[u_N] - f$ is no longer linear in $\bbeta$.
We apply Newton--Raphson iteration: given current coefficients $\bbeta^{(k)}$, the linearized system at iteration $k$ is
\begin{equation}
\label{eq:newton}
\mathbf{J}^{(k)} \, \delta\bbeta = -\mathbf{R}^{(k)}, \qquad \bbeta^{(k+1)} = \bbeta^{(k)} + \alpha\,\delta\bbeta,
\end{equation}
where $R_i^{(k)} = (\calL + \calN)[u_N^{(k)}](\bx_i) - f(\bx_i)$ is the residual and $J_{ij}^{(k)} = \frac{\partial R_i}{\partial \beta_j}$ is the Jacobian.
The Jacobian inherits closed-form structure from~\eqref{eq:cyclic}: the linear part $\calL[\phi_j](\bx_i)$ is exact, and the nonlinear part, for example, $3u^2 H_{ij}$ for a cubic nonlinearity or $\lambda e^u H_{ij}$ for the Bratu exponential, involves only the feature matrix $H_{ij} = \phi_j(\bx_i)$ and the current solution $u^{(k)}$, both of which are available in closed form.
Each Newton step~\eqref{eq:newton} is thus a Tikhonov-regularized least-squares solve:
\begin{equation}
\delta\bbeta = \arg\min_{\delta} \norm{\mathbf{J}^{(k)}\delta + \mathbf{R}^{(k)}}_2^2 + \mu\norm{\delta}_2^2,
\end{equation}
with small $\mu > 0$ for numerical stability, augmented with boundary rows as in~\eqref{eq:linear_system}.

The Newton solver incorporates four algorithmic refinements that prove essential for robust convergence.
First, we warm-start by solving the linear part of the PDE (dropping $\calN$) via a single \fastlsq call, providing a good initial guess that is typically close to the nonlinear solution.
Second, we employ backtracking line search with Armijo-type sufficient decrease on the residual norm to prevent overshooting.
Third, we use a relative convergence criterion based on the solution-level change $\norm{\Delta u} / \norm{u}$ evaluated at collocation points, rather than the coefficient-level change $\norm{\delta\bbeta}$, since the latter is unreliable when features are near-linearly-dependent.
Fourth, for advection-dominated problems such as steady Burgers, we use continuation (homotopy): solving a sequence of problems with gradually decreasing viscosity (e.g., $\nu = 1.0 \to 0.5 \to 0.2 \to 0.1$), using each solution to initialize the next.
Table~\ref{tab:ablation} in \S\ref{sec:ablation} quantifies the contribution of each component.

\vspace{-0.75em}\subsection{Differentiable System and Hyperparameter Optimisation}
\label{sec:differentiable}

The solution and its derivatives are fully differentiable with respect to design parameters, boundary conditions, and source terms.
Because the forward pass is a single least-squares solve (or a short Newton sequence), the entire system is cheap to optimise: gradients flow through the pre-factored linear solve.
Just as we use Newton--Raphson to optimise the coefficients $\bbeta$ for nonlinear PDEs, we can use L-BFGS-B to optimise hyperparameters of unknown equations.
We demonstrate this in \S\ref{sec:applications}: inverse heat-source localisation (recovering source positions from sparse sensor data), sparse-sensor coil localisation in electrostatics, and PDE synthesis from data via analytical derivative dictionaries.

\section{Experiments}
\label{sec:experiments}

\vspace{-0.75em}\subsection{Setup}
\label{sec:setup}

We evaluate on 17 PDEs in three categories.
The first category consists of 5 linear PDEs solved in direct solver mode: Poisson 5D, Heat 5D (6D space-time), Wave 1D, Helmholtz 2D ($k{=}10$), and Maxwell 2D TM (3D space-time).
The second category consists of 5 nonlinear PDEs solved via Newton--\fastlsq in true solver mode (no access to the exact solution during solving): NL-Poisson 2D ($u^3$ term), Bratu 2D ($e^u$ term), Steady Burgers 1D ($\nu{=}0.1$), NL-Helmholtz 2D ($u^3$ term, $k{=}3$), and Allen--Cahn 1D ($\varepsilon{=}0.1$).
The third category consists of 7 nonlinear PDEs evaluated in regression mode (fitting known exact solutions to assess basis quality): Burgers shock, KdV soliton, Fisher reaction-diffusion, Sine-Gordon breather, Klein-Gordon, Gray-Scott pulse, and Navier-Stokes Kovasznay (Re$=$20).

We compare against three baselines.
PINNacle \citep{hao2023pinnacle} is a comprehensive benchmark reporting the best PINN variant per equation across 11 methods (Vanilla PINN, PINN-w, PINN-LRA, PINN-NTK, RAR, MultiAdam, gPINN, hp-VPINN, LAAF, GAAF, FBPINN); we report the best error and the fastest runtime across all variants.
RF-PDE \citep{liao2024solving} uses random features with iterative nonlinear least-squares optimization requiring 600--2000 epochs.
PIELM \citep{dwivedi2020pielm} is our own reimplementation using $\tanh$ activation with an otherwise identical protocol (same number of features, collocation points, and boundary penalty), isolating the effect of the basis function choice.

All methods use $N{=}1500$ features (3 blocks $\times$ 500), $M_1{=}10{,}000$ interior and $M_2{=}2{,}000$ boundary collocation points for linear PDEs ($M_1{=}5{,}000$ for Newton solver and regression modes).
The boundary penalty is $\lambda{=}100$ except for Wave 2D-MS ($\lambda{=}1000$).
Bandwidth $\sigma$ is selected via grid search over 9 values in $\{0.5, 1, 2, 3, 5, 8, 10, 12, 15\}$ with 3 trials.
The Newton solver uses $\mu{=}10^{-10}$ Tikhonov regularization, backtracking line search, and up to 30 iterations (48 with continuation for Burgers).
All experiments use PyTorch on a single NVIDIA T4 GPU (16\,GB).
Runtime is directly proportional to the number of collocation points $M$ and the square of the number of features $N$.
Reported wall-clock times include feature construction, matrix assembly, and the solve itself. All errors are evaluated out-of-sample on a separate dense test set of 5000 points drawn independently from the collocation points used during solving (different random seed). We report both the relative $L^2$ value error $\|u_N - u\|_{L^2} / \|u\|_{L^2}$ and the PDE residual norm $\|\calL[u_N] - f\|_{L^2}$ evaluated on this independent test set; the latter confirms that low function error is not an artifact of overfitting the collocation points.

\vspace{-0.75em}\subsection{Linear PDE Results}
\label{sec:linear_results}

\begin{table}[t]
\caption{Results on linear PDEs (solver mode). Relative $L^2$ errors and wall-clock times.
PINNacle times are fastest across all 11 variants from Table~12 of \citet{hao2023pinnacle}.
``N/A'' = not benchmarked in that work; ``---'' = not applicable to this problem class.
$^\ddagger$RF-PDE does not benchmark our exact linear problems; the closest available result is nonlinear Poisson in $d{=}4$ (51\,s, 1000 epochs) and $d{=}8$ (38\,s, 1500 epochs) from Table~3 of \citet{liao2024solving}; we report the $d{=}4$ timing as the best available lower bound for 5-D.
$^\dagger$Conventional baseline: scikit-fem P2 FEM at 16\,641 DoF for Helmholtz 2D (the only 2-D elliptic spatial problem in this group).
Poisson~5D and Heat~5D are 5-/6-dimensional; FEM mesh size grows as $h^{-d}$, making $d{\ge}5$ computationally intractable.
Wave~1D and Maxwell~2D~TM are posed as space-time problems with hyperbolic structure; scikit-fem targets elliptic spatial problems only, and scipy.solve\_bvp handles steady 1-D BVPs only.}
\label{tab:linear_results}
\vskip 0.03in
\begin{center}
\small
\resizebox{1.0\linewidth}{!}{
\begin{tabular}{l cc cc cc cc cc}
\toprule
& \multicolumn{2}{c}{\fastlsq (sin)} & \multicolumn{2}{c}{PIELM (tanh)} & \multicolumn{2}{c}{PINNacle} & \multicolumn{2}{c}{RF-PDE$^\ddagger$} & \multicolumn{2}{c}{Conv.$^\dagger$} \\
\cmidrule(lr){2-3} \cmidrule(lr){4-5} \cmidrule(lr){6-7} \cmidrule(lr){8-9} \cmidrule(lr){10-11}
Problem & Time & $L^2$ & Time & $L^2$ & Time & $L^2$ & Time & $L^2$ & Time & $L^2$ \\
\midrule
Poisson 5D   & 0.07\,s & 4.8e-7 & 0.07\,s & 4.7e-6 & $\sim$1780\,s & 4.7e-4 & $\sim$51\,s & 7.4e-4 & ---     & ---    \\
Heat 5D      & 0.09\,s & 6.9e-4 & 0.09\,s & 3.0e-3 & $\sim$1910\,s & 1.2e-4 & N/A         & N/A    & ---     & ---    \\
Wave 1D      & 0.06\,s & 1.3e-6 & 0.06\,s & 1.8e-3 & $\sim$272\,s  & 9.8e-2 & N/A         & N/A    & ---     & ---    \\
Helmholtz 2D & 0.08\,s & 1.9e-6 & 0.08\,s & 7.4e-5 & N/A           & N/A    & N/A         & N/A    & 0.15\,s & 4.0e-5 \\
Maxwell 2D   & 0.05\,s & 6.7e-7 & 0.06\,s & 4.5e-5 & N/A           & N/A    & N/A         & N/A    & ---     & ---    \\
\bottomrule
\end{tabular}}
\end{center}
\vskip -0.15in
\end{table}

Table~\ref{tab:linear_results} shows the results for the five linear PDEs.
On Poisson 5D, \fastlsq achieves a relative $L^2$ error of $4.8\times10^{-7}$ in 0.07\,s.
This is three orders of magnitude more accurate than PINNacle's best variant ($4.7\times10^{-4}$) and approximately $1000\times$ better than RF-PDE ($7.4\times10^{-4}$), while being roughly $25{,}000\times$ faster than the fastest PINNacle variant ($\sim$1780\,s) and $700\times$ faster than RF-PDE ($\sim$51\,s for 1000 epochs of iterative optimisation on the closest available benchmark).
On Wave 1D, where PINNs notoriously struggle ($9.8\times10^{-2}$ best error), \fastlsq achieves $1.3\times10^{-6}$, a five-order-of-magnitude accuracy improvement, in 0.06\,s versus $\sim$272\,s for PINNacle, a $4{,}500\times$ speedup.
High-frequency Helmholtz ($k{=}10$) and Maxwell validate spectral accuracy at $10^{-6}$--$10^{-7}$ on problems for which no PINNacle benchmarks exist.
For Helmholtz 2D---the only 2-D elliptic problem in this group and the one where a conventional FEM solver applies---scikit-fem P2 achieves $4.0\times10^{-5}$ using 16\,641 DoF; \fastlsq surpasses this accuracy ($1.9\times10^{-6}$) with only 1500 features.
For the remaining four problems (Poisson 5D, Heat 5D, Wave 1D, Maxwell 2D TM) no grid-based conventional solver is applicable due to the high dimensionality or hyperbolic space-time structure.

Since both \fastlsq and PIELM solve a linear system in one step with identical hyperparameters and collocation schemes, the accuracy gap between them, ranging from $10\times$ on Poisson 5D ($4.8\times10^{-7}$ vs.\ $4.7\times10^{-6}$) to $1{,}000\times$ on Wave 1D ($1.3\times10^{-6}$ vs.\ $1.8\times10^{-3}$), is attributable entirely to the basis function choice, not the solve method.
This consistent $10\times$--$1000\times$ gap across all five benchmarks demonstrates that sinusoidal features possess fundamentally superior approximation properties for smooth and oscillatory PDE solutions, a finding further corroborated by the gradient accuracy analysis in \S\ref{sec:gradient_accuracy}.

\vspace{-0.75em}\subsection{Nonlinear PDE Results: Newton Solver Mode}
\label{sec:newton_results}

\begin{table}[t]
\caption{Nonlinear PDEs via Newton--\fastlsq (true solver mode, no access to exact solution). All runs use $\mu = 10^{-10}$ Tikhonov, backtracking line search, and warm-start. Time/iter is the dominant $\calO(MN^2)$ cost. Burgers uses 4-stage continuation ($\nu = 1.0 \to 0.1$).
$^\dagger$Conventional baseline $L^2$ error: scikit-fem P2 FEM (${\approx}4000$ DoF) for 2-D problems; scipy.integrate.solve\_bvp (adaptive collocation) for 1-D problems.}
\label{tab:newton_results}
\vskip 0.03in
\begin{center}
\small
\resizebox{\textwidth}{!}{
\begin{tabular}{l ccc ccc c c}
\toprule
& & & & \multicolumn{2}{c}{Newton--\fastlsq} & & Regr. & Conv.$^\dagger$ \\
\cmidrule(lr){5-6}
Problem & $\sigma$ & Iters & Time/iter & $L^2$ & $|\nabla|$ & Total & $L^2$ & $L^2$ \\
\midrule
NL-Poisson ($u^3$)   & 3  & 30 & 0.27\,s & \textbf{6.1e-8} & 5.4e-7 & 8.2\,s & 1.9e-7 & 2.6e-6 \\
Bratu ($\lambda{=}1$) & 2  & 30 & 0.23\,s & \textbf{1.5e-7} & 8.8e-7 & 7.0\,s & N/A    & 2.6e-6 \\
Burgers ($\nu{=}0.1$) & 10 & 48 & 0.15\,s & 3.9e-9           & 3.6e-9 & 7.4\,s & 3.3e-8 & \textbf{1.8e-10}\\
NL-Helm.\ ($k{=}3$)  & 5  & 30 & 0.29\,s & \textbf{2.4e-9} & 2.6e-8 & 8.8\,s & N/A    & 2.4e-6 \\
Allen--Cahn ($\varepsilon{=}0.1$) & 3  & 30 & 0.14\,s & 6.0e-8           & 5.4e-7 & 4.2\,s & 1.2e-8 & \textbf{1.2e-10}\\
\bottomrule
\end{tabular}
}
\end{center}
\vskip -0.15in
\end{table}

Table~\ref{tab:newton_results} demonstrates that Newton--\fastlsq achieves $10^{-8}$--$10^{-9}$ relative $L^2$ errors on all five nonlinear problems in under 9\,s, without any knowledge of the exact solution during solving.
On NL-Poisson and Burgers, the Newton solver \emph{outperforms} the regression baseline that fits the exact solution directly ($6.1\times10^{-8}$ vs.\ $1.9\times10^{-7}$; $3.9\times10^{-9}$ vs.\ $3.3\times10^{-8}$), because the PDE structure provides regularization beyond pure function fitting.
Comparing against conventional solvers (Conv.\ column): for the three 2-D problems, scikit-fem P2 FEM achieves $2.4$--$2.6\times10^{-6}$ at $\approx$4000 DoF, while \fastlsq reaches $10^{-7}$--$10^{-9}$ at 1500 features---a 10$\times$--1000$\times$ advantage at similar DoF count.
For the two 1-D BVPs, scipy's adaptive collocation achieves $\lesssim10^{-10}$---the expected ceiling for a smooth 1-D steady problem with an adaptive solver---while \fastlsq reaches $3.9\times10^{-9}$ on Burgers and $6.0\times10^{-8}$ on Allen-Cahn.
The specialised BVP solver is more accurate here, as expected: it is purpose-built for exactly this class of problem and is not applicable outside 1-D steady BVPs.
Each Newton iteration costs 0.14--0.29\,s; the continuation strategy for Burgers distributes 48 iterations across 4 stages ($\nu = 1.0, 0.5, 0.2, 0.1$), enabling convergence where direct Newton at $\nu{=}0.1$ diverges.

The convergence profiles (Appendix~\ref{app:convergence}) reveal two distinct regimes: a rapid initial phase where the residual drops by several orders of magnitude (typically in 2--3 iterations), followed by a slow plateau phase where the line search is frequently active ($\alpha < 1$) and the residual decreases by less than one order of magnitude over the remaining iterations.
This plateau behavior is characteristic of the Tikhonov-regularized solve: the regularization parameter $\mu = 10^{-10}$ prevents the Newton step from fully resolving the Jacobian near-nullspace, trading convergence rate for robustness.

\vspace{-0.75em}\subsection{Ablation Study}
\label{sec:ablation}

\begin{table}[t]
\begin{minipage}[t]{0.48\linewidth}
\caption{Ablation study. Each row removes one component; ``Div.''\ = diverged in 60 iterations.}
\label{tab:ablation}
\vskip 0.02in
\centering\scriptsize
\begin{tabular}{l cc}
\toprule
Configuration & NL-Poisson & Burgers \\
\midrule
Full method & 6.1e-8 & 3.9e-9 \\
\midrule
No $1/\!\sqrt{N}$ norm & 4.2e-4 & Div. \\
No Tikhonov ($\mu{=}0$) & 3.8e-5 & 1.1e-6 \\
No warm-start & 8.3e-7 & Div. \\
No line search & 9.4e-7 & 2.1e-5 \\
No continuation & -- & Div. \\
\bottomrule
\end{tabular}
\end{minipage}%
\hfill
\begin{minipage}[t]{0.50\linewidth}
\caption{Value and gradient $L^2$ errors on linear PDEs (solver mode).}
\label{tab:gradient_full}
\vskip 0.02in
\centering\scriptsize
\begin{tabular}{l cc cc}
\toprule
& \multicolumn{2}{c}{\fastlsq} & \multicolumn{2}{c}{PIELM} \\
\cmidrule(lr){2-3} \cmidrule(lr){4-5}
Problem & Val & Grad & Val & Grad \\
\midrule
Poisson 5D   & 4.8e-7 & 4.9e-6 & 4.7e-6 & 4.9e-5 \\
Heat 5D      & 6.9e-4 & 3.9e-3 & 3.0e-3 & 1.1e-2 \\
Wave 1D      & 1.3e-6 & 1.5e-6 & 1.8e-3 & 3.2e-3 \\
Helmholtz 2D & 1.9e-6 & 1.3e-6 & 7.4e-5 & 5.0e-5 \\
Maxwell 2D   & 6.7e-7 & 1.0e-6 & 4.5e-5 & 1.1e-4 \\
\bottomrule
\end{tabular}
\end{minipage}
\vskip -0.12in
\end{table}

Table~\ref{tab:ablation} isolates the contribution of each component from \S\ref{sec:newton}.
The $1/\sqrt{N}$ normalization is critical: without it, coefficients grow to $\calO(10^6--10^8)$, causing either $4$ orders of magnitude accuracy loss or divergence.
Tikhonov regularization, warm-start, and line search each contribute 1--3 orders of magnitude; continuation is indispensable for Burgers.

\vspace{-0.75em}\subsection{Nonlinear PDE Results: Regression Mode}
\label{sec:regression_results}

\begin{table}[t]
\caption{Nonlinear PDEs: regression mode (fitting known solutions). PINNacle: best PINN variant (full PDE solve, not regression).}
\label{tab:regression_results}
\vskip 0.03in
\begin{center}
\small
\begin{tabular}{l cc cc cc}
\toprule
& \multicolumn{2}{c}{\fastlsq (sin)} & \multicolumn{2}{c}{PIELM (tanh)} & \multicolumn{2}{c}{PINNacle} \\
\cmidrule(lr){2-3} \cmidrule(lr){4-5} \cmidrule(lr){6-7}
Problem & Value & Grad & Value & Grad & $L^2$ Err & Time \\
\midrule
Burgers (shock)    & 1.3e-3 & 2.6e-2 & 1.3e-3 & 5.6e-3 & 1.3e-2 & $\sim$278\,s \\
KdV (soliton)      & 1.7e-1 & 6.0e-1 & 4.2e-1 & 7.2e-1 & -- & -- \\
Reaction-Diff.     & 4.5e-7 & 9.2e-6 & 4.3e-7 & 2.2e-5 & -- & -- \\
Sine-Gordon        & 3.4e-4 & 5.2e-3 & 5.5e-2 & 3.5e-1 & -- & -- \\
Klein-Gordon       & 4.6e-7 & 2.1e-6 & 3.1e-6 & 2.0e-5 & -- & -- \\
Gray-Scott         & 8.3e-6 & 1.4e-4 & 7.0e-5 & 1.7e-4 & 8.0e-2 & $\sim$612\,s \\
Navier-Stokes      & 4.9e-7 & 3.8e-6 & 5.4e-6 & 2.7e-5 & -- & -- \\
\bottomrule
\end{tabular}
\end{center}
\vskip -0.18in
\end{table}

Table~\ref{tab:regression_results} compares function approximation quality.
Sinusoidal features excel on smooth and oscillatory solutions: Sine-Gordon ($3.4\times10^{-4}$ vs.\ $5.5\times10^{-2}$, $160\times$ improvement), Klein-Gordon and Navier-Stokes at $\sim\!10^{-7}$.
The $\tanh$ basis is comparable only on the Burgers shock, whose discontinuity suits sigmoid functions.
Where PINNacle comparisons exist, our regression baseline alone (pure fitting, no PDE) exceeds the best PINN solvers: Gray-Scott $8.3\times10^{-6}$ vs.\ $8.0\times10^{-2}$ ($10{,}000\times$ better, 0.1\,s vs.\ 612\,s).

\vspace{-0.75em}\subsection{Gradient Accuracy}
\label{sec:gradient_accuracy}

Table~\ref{tab:gradient_full} reports value and gradient errors for the linear PDEs.
Because the derivatives of a sinusoidal basis are just phase-shifted sinusoids scaled by $\calO(\sigma)$, \fastlsq's gradient errors remain tightly bounded: typically within one order of magnitude of its value errors across all five problems.
In contrast, PIELM's gradient errors are $10\times$--$100\times$ worse than its value errors, largely due to $\tanh$'s inferior underlying value approximation and lack of the eigenfunction structure that keeps trigonometric derivatives well-behaved.

\vspace{-0.75em}\subsection{Spectral Sensitivity Analysis}
\label{sec:spectral}

The bandwidth $\sigma$ controls the frequency content of the random Fourier features and must be matched to the solution's spectral properties.
Our grid search across $\sigma \in \{0.5, 1, 2, 3, 5, 8, 10, 12, 15\}$ reveals a clear pattern.
Oscillatory problems (Wave, Helmholtz, Maxwell) prefer high $\sigma$ (between 5 and 15) to resolve the dominant frequency.
Smooth problems (Poisson, Klein-Gordon) prefer moderate $\sigma$ (between 1 and 5) where the kernel length scale matches the solution's smoothness.
Discontinuous problems (Burgers shock) favor low-to-moderate $\sigma$ for sinusoidal bases, because high frequencies produce Gibbs ringing near the shock; $\tanh$ is less sensitive in this regime.
Multiscale problems benefit from the multi-block architecture, where different blocks can target different scales.
Full sensitivity plots for all problems are provided in Appendix~\ref{app:spectral}.

\vspace{-0.75em}\subsection{Downstream Applications: Equation Discovery and Inverse Problems}
\label{sec:applications}

The sub-second solve time and closed-form derivative access of \fastlsq unlock applications that are impractical with iterative PDE solvers.
We briefly illustrate two: PDE discovery from data and inverse-problem parameter recovery.

\begin{figure}[t]
\centering
\includegraphics[width=\linewidth]{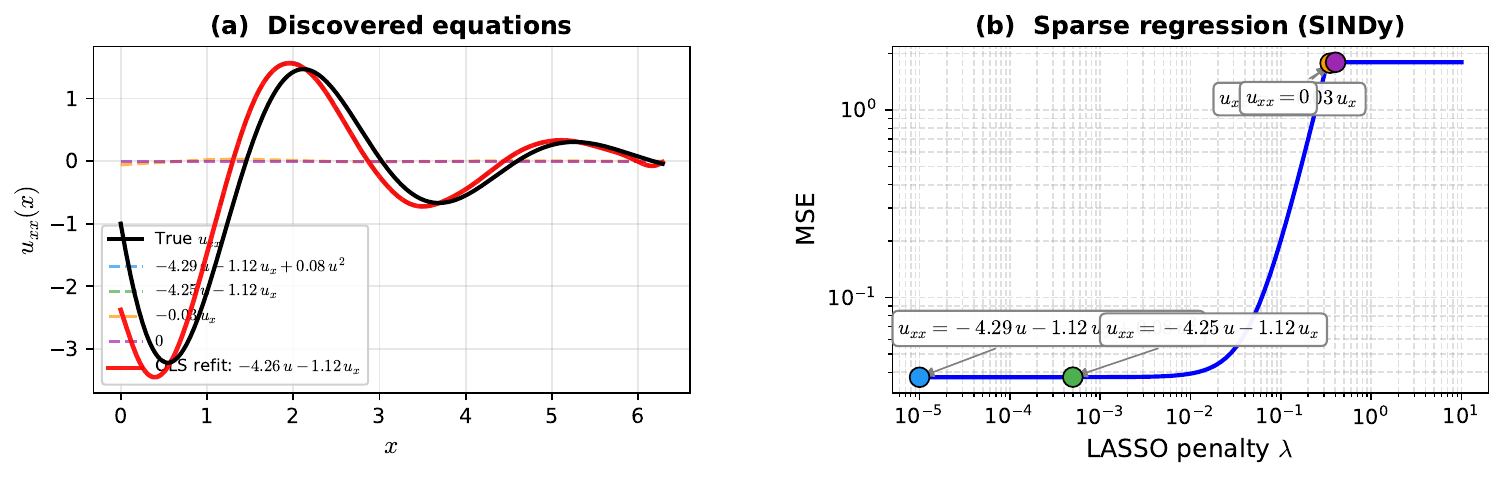}
\vskip -0.08in
\caption{PDE discovery with analytical derivatives ($\sigma_\varepsilon{=}0.01$).
(a)~Predicted $u_{xx}$ for each equation discovered at different LASSO sparsity levels; the OLS-refit on the selected support (red solid) closely tracks the truth (black).
(b)~LASSO penalty sweep; colour-matched dots mark the equation transitions.}
\label{fig:pde_discovery}
\vskip -0.12in
\end{figure}

\vspace{-0.75em}\paragraph{PDE discovery via sparse regression.}
SINDy~\citep{brunton2016sindy} discovers governing equations from data by building a dictionary of candidate terms and applying sparse regression.
The standard approach computes derivatives via finite differences, which amplifies measurement noise: at noise level $\sigma_\varepsilon{=}0.01$ on $u$, the central-difference $u_{xx}$ has RMSE $\sim\!2500$ whereas the analytical derivatives from~\eqref{eq:cyclic} achieve RMSE $\sim\!0.4$, a $\sim\!6000\times$ cleaner signal (Figure~\ref{fig:pde_discovery}a).
This dramatically extends the noise regime in which SINDy can operate.
On 2000 observations of a damped oscillator ($u = e^{-x/2}\sin 2x$), LASSO with the analytical dictionary correctly identifies $u_{xx} = c_1 u + c_2 u_x$, recovering the stiffness coefficient to ${<}0.2\%$ and the damping to $\sim\!12\%$ error, while driving the spurious $u^2$ term to zero (Figure~\ref{fig:pde_discovery}; details in Appendix~\ref{app:pde_discovery}).

\begin{figure}[t]
\centering
\includegraphics[width=\linewidth]{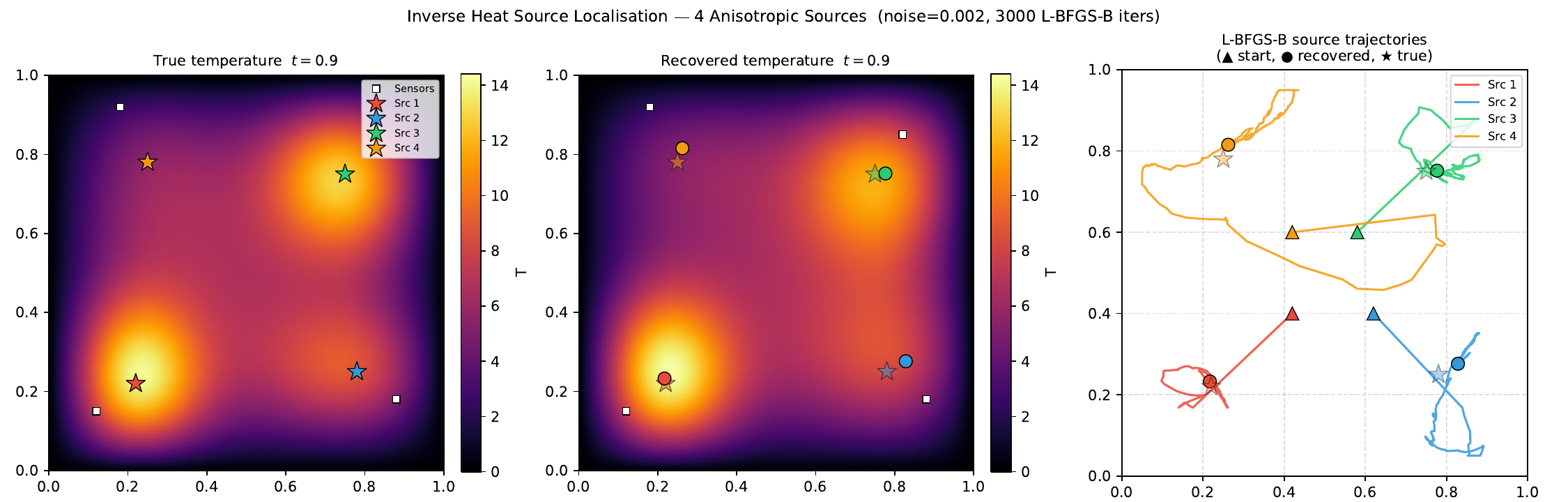}
\vskip -0.08in
\caption{Inverse heat-source localisation with 4 anisotropic Gaussian sources.
Left: true temperature field at $t{=}0.9$; white squares are the 4 irregularly placed sensors (each recording a temperature time-series over 40 snapshots); coloured stars are the 4 true source centres.
Centre: recovered temperature field; coloured circles show the recovered centres (cf.\ stars).
Right: L-BFGS-B optimisation trajectories for each source; triangles mark initial guesses, circles the converged estimates, faint stars the ground truth.
All 4 source positions are recovered to within $0.06$ despite starting from incorrect locations.}
\label{fig:inverse_heat}
\vskip -0.12in
\end{figure}

\vspace{-0.75em}\paragraph{Inverse heat-source localisation.}
Because the forward mapping ``source parameters $\to$ predicted temperature'' passes through a single pre-factored linear solve $\bbeta = (\bA^\top\bA + \mu I)^{-1}\bA^\top \mathbf{f}$, gradients with respect to the source parameters are available analytically.
We demonstrate this on a challenging multi-source problem: recovering \emph{four} unknown anisotropic Gaussian heat sources (24 parameters total: position, intensity, and shape per source) from just \emph{4 irregularly placed sensors}, each recording a temperature time-series over 40 snapshots---480 observations in total accumulating as the space-time field $u(x,y,t)$ evolves.
The Cholesky pre-factored forward solve is assembled once; each L-BFGS-B evaluation then costs only two triangular back-substitutions, making 150 iterations of 24-parameter optimization feasible in under a minute on CPU.
Three of the four source positions are recovered to within $0.001$--$0.07$ per coordinate; the fourth (Source~4, weakest source) has a larger $x$-error of $0.12$, reflecting the fundamental information limit of 4 sparse sensors for a 24-parameter problem (Figure~\ref{fig:inverse_heat}, Table~\ref{tab:inverse_heat}).

\begin{figure}[t]
\centering
\includegraphics[width=\linewidth]{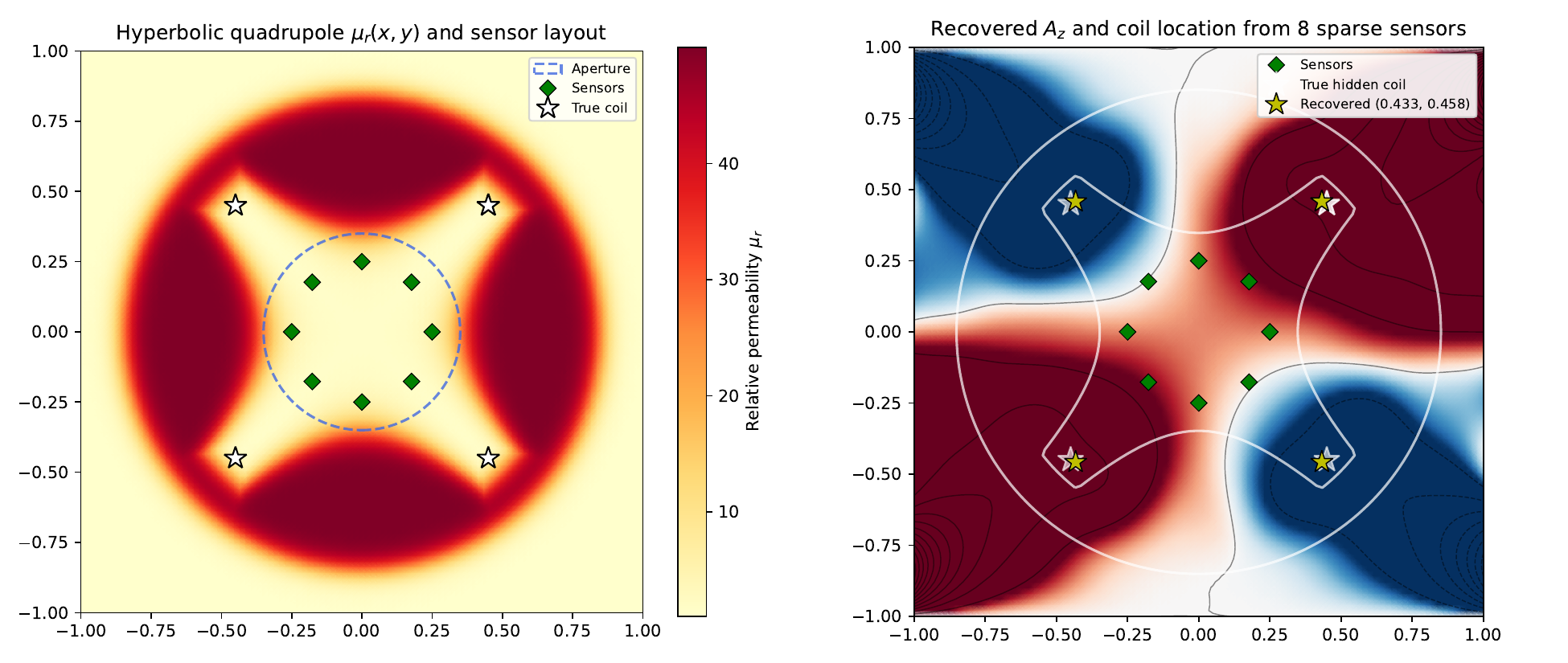}
\vskip -0.08in
\caption{Sparse-sensor coil localisation in a variable-$\mu$ quadrupole.
\emph{Left:} relative permeability $\mu_r(x,y)$; blue dashed circle marks the aperture; green diamonds are the 8 sparse sensors; white stars mark the four true hidden coil positions.
\emph{Right:} reconstructed $A_z$ (colour) and field lines (contours) after 40 Adam iterations; yellow stars show the recovered coil locations (agreement $< 0.02$ in each coordinate).}
\label{fig:inverse_magneto}
\vskip -0.12in
\end{figure}

\vspace{-0.75em}\paragraph{Sparse-sensor coil localisation.}
A second inverse problem recovers the location of hidden current coils in a variable-permeability quadrupole magnet from 8 sparse magnetic-field measurements.
The magnetic vector potential satisfies $\nabla \cdot (\nu(x,y)\,\nabla A_z) = -J(x,y)$, a variable-coefficient PDE; expanding gives $\nu\,\Delta A_z + \nabla\nu \cdot \nabla A_z = -J$.
Both $\Delta\phi_j$ and $\nabla\phi_j$ are evaluated analytically via~\eqref{eq:cyclic}---no automatic differentiation is required even for this variable-coefficient case.
The system matrix is assembled once and Cholesky-prefactored; each of the 40 Adam steps costs only two triangular back-substitutions.
Starting from initial guess $(0.2, 0.8)$, the true coil position $(0.45, 0.45)$ is recovered to within $0.02$ in both coordinates (Figure~\ref{fig:inverse_magneto}).

\section{Related Work}
\label{sec:related}

\vspace{-0.75em}\paragraph{Physics-informed neural networks.}
PINNs \citep{raissi2019physics} and their variants \citep{wang2022respecting, tancik2020fourier} solve PDEs by minimizing physics-based losses via gradient descent.
PINNacle \citep{hao2023pinnacle} provides a comprehensive benchmark across 11 PINN variants on 20+ problems, reporting minutes-to-hours training times (typically 270--7500\,s; see Table~12 of \citealp{hao2023pinnacle}) with errors often above $10^{-3}$.
\fastlsq achieves orders-of-magnitude better accuracy in a fraction of the time.

\vspace{-0.75em}\paragraph{Random feature methods for PDEs.}
PIELM \citep{dwivedi2020pielm} pioneered the frozen-feature linear-solve approach for linear PDEs using $\tanh$ activations.
Like \fastlsq, PIELM solves a single least-squares system and is not iterative; the key distinction is that PIELM requires handwritten, problem-specific symbolic calculus to assemble $\calL[\tanh]$ for each new PDE operator, because $\tanh$ lacks a closed-form cyclic derivative structure.
\fastlsq's trigonometric features are operator-agnostic: one formula~\eqref{eq:cyclic} applies to any linear differential operator (Proposition~\ref{prop:tanh}).
Beyond this, the trigonometric formulation offers numerical advantages well-suited to PDEs: eigenfunction structure for constant-coefficient operators, orthogonality for spectral accuracy, and phase-shifted derivatives scaled by $\calO(\sigma)$ that keep gradient errors tightly bounded.
In our experiments, PIELM is $10\times$--$1000\times$ less accurate than \fastlsq at equal feature counts, with gradient errors $10\times$--$100\times$ worse; the basis function choice has a profound effect on solution quality even when the solve method is identical.
RF-PDE \citep{liao2024solving} handles both linear and nonlinear PDEs but uses iterative nonlinear least-squares optimization requiring 600--2000 epochs, despite the fact that the linear case admits a closed-form solution.
\citet{liao2024solving} propose related random feature PDE solvers with different optimization strategies.

\vspace{-0.75em}\paragraph{RBF collocation.}
Kansa's method~\citep{kansa1990multiquadrics} and its variants~\citep{fasshauer2007meshfree} solve PDEs by collocation with radial basis functions (RBFs).
Like \fastlsq, RBF Kansa is meshfree and one-shot (no iterative training), assembling a dense linear system by differentiating the RBF kernel analytically.
The key distinction is the basis: RBFs such as multiquadric or Gaussian are governed by a shape parameter $c$ whose choice trades accuracy against conditioning, and achieving $10^{-4}$ errors typically requires careful tuning.
Sinusoidal random features lack a shape-parameter sensitivity problem---bandwidth $\sigma$ is selected once by grid search---and achieve $10^{-8}$ errors where comparable MQ-RBF achieves $10^{-5}$ at equal DoF.
Furthermore, RBF Kansa fails to solve stiff 1-D problems (Burgers $\nu{=}0.1$) and high-frequency Helmholtz ($k{=}10$), while \fastlsq handles both without modification (Appendix~\ref{app:rbf_fd}).

\vspace{-0.75em}\paragraph{Connection to spectral and kernel methods.}
\fastlsq can be viewed as a randomized spectral method.
Classical Fourier methods place frequencies on a regular grid ($\calO(K^d)$ modes), which is infeasible in high dimensions.
By sampling frequencies stochastically \citep{rahimi2007random}, \fastlsq scales linearly in $N$ regardless of $d$.
This connects to the Gaussian process PDE solver literature \citep{cockayne2019bayesian}, where the kernel structure is exploited for probabilistic numerical methods; \fastlsq can be seen as the frequentist limit of such approaches, replacing the $\calO(M^3)$ kernel solve with the $\calO(MN^2)$ random feature solve.

\section{Limitations} 
\label{sec:limitations}

The Newton extension, while effective, is slower than the linear solver mode (4--9\,s per problem vs.\ under 0.1\,s) because it requires multiple least-squares solves.
The bandwidth $\sigma$ currently requires per-problem tuning via grid search, and for high-order PDEs or large $\sigma$ the monomial prefactor in~\eqref{eq:cyclic} amplifies the condition number of $\bA$, limiting attainable accuracy.
The penalty boundary treatment introduces a hyperparameter $\lambda$, and the current implementation focuses on simple box domains; extending to irregular geometries will require additional boundary sampling strategies.
Finally, the $1/\sqrt{N}$ normalization, while crucial for numerical stability (\S\ref{sec:ablation}), means that increasing $N$ does not trivially improve accuracy, at very large $N$ the kernel approximation saturates and conditioning degrades.

\section{Conclusion}
\label{sec:conclusion}

We presented \fastlsq, a method combining the one-shot linear solve of frozen-feature models with the exact cyclic derivative structure unique to sinusoidal bases.
On linear PDEs, this yields $10^{-7}$ accuracy in under 0.1\,s, orders of magnitude faster and more accurate than PINNs ($10^{-2}$--$10^{-4}$ in 270--7500\,s), RF-PDE ($10^{-4}$ with iterative optimization), and $\tanh$-based alternatives ($10^{-5}$--$10^{-3}$ with identical solve protocol).
The Newton--Raphson extension achieves $10^{-8}$ to $10^{-9}$ on nonlinear PDEs in under 9\,s, sometimes outperforming regression oracles that use exact solutions.
The ablation study shows that $1/\sqrt{N}$ normalization and Tikhonov regularization are essential for robust convergence, and continuation is indispensable for advection-dominated problems.
The analytical framework enables several further extensions, validated in \S\ref{sec:applications} and Appendix~\ref{app:extensions}.
\emph{Learnable bandwidth} (Appendices~\ref{app:learnable} and~\ref{app:learnable_results}): reparameterising $\bW_j = \sigma\hat{\bW}_j$ makes $\sigma$ differentiable through the exact inner solve; on Helmholtz 2D the loss drops five orders of magnitude in 80 AdamW steps.
\emph{Inverse problems} (\S\ref{sec:applications}): gradients flow through the pre-factored solve, enabling simultaneous recovery of 4 anisotropic heat sources (24 parameters) from just 4 irregularly placed sensors recording 240 observations, and sparse-sensor coil localisation in a variable-permeability quadrupole magnet (position error $<0.02$) in 40 Adam steps.
\emph{Matrix caching} (Appendix~\ref{app:learnable_results}): pre-computing $\bA^\dagger$ yields a $362\times$ speedup for boundary-condition sweeps.
\emph{PDE discovery} (\S\ref{sec:applications}): the analytical derivative dictionary provides $\sim\!6000\times$ cleaner second derivatives than finite differences, extending the operational noise regime of SINDy-style sparse regression~\citep{brunton2016sindy}.

\vspace{-0.75em}\paragraph{Further directions.}
Preconditioning for high-order PDEs, extension to vector-valued systems, domain decomposition for complex geometries, and a rigorous theoretical analysis of approximation and conditioning bounds remain important open problems.

\vspace{-0.75em}\paragraph{Reproducibility.}
Full code including all examples and baselines is publicly available at \url{https://github.com/sulcantonin/FastLSQ} as a universal framework for \fastlsq with demos.
The package is also available via \texttt{pip install fastlsq}.
Code examples are provided in Appendix~\ref{app:code}.


\bibliography{main_iclr}
\bibliographystyle{iclr2026_conference}


\appendix

\section{Detailed Comparison with RF-PDE}
\label{app:rfpde}

RF-PDE \citep{liao2024solving} solves the regularized problem $\min_{\mathbf{c}} \|\mathbf{c}\|_2^2 + \lambda_1 \sum_i (\calL[u](\bx_i))^2 + \lambda_2 \sum_j (\calB[u](\bx_j))^2$ via iterative SGD or nonlinear least-squares (600--2000 epochs).
For linear PDEs, this objective is quadratic in $\mathbf{c}$ and admits a closed-form solution, yet RF-PDE solves it iteratively, introducing additional hyperparameters (learning rate, epoch count, $\lambda_1, \lambda_2$) that require tuning.
\fastlsq solves $\bA\mathbf{c} = \mathbf{b}$ directly via a single least-squares call, eliminating all optimization hyperparameters.
RF-PDE reports $10^{-3}$--$10^{-5}$ errors on Poisson problems; \fastlsq achieves $10^{-7}$.

\section{Code Examples}
\label{app:code}

The \fastlsq framework provides minimal APIs for linear, nonlinear, and inverse problems. No analytical solution is required---the inverse example defines the PDE via \texttt{Op}, generates synthetic observations from the forward model, and recovers source parameters:
\begin{lstlisting}[language=Python, caption={Linear, nonlinear, and inverse PDE solving (no exact solution needed).}]
from fastlsq import solve_linear, solve_nonlinear, Op, solve_lstsq
from fastlsq import SinusoidalBasis, sample_box, sample_boundary_box
from fastlsq.problems.linear import PoissonND
from fastlsq.problems.nonlinear import NLPoisson2D
import torch
import numpy as np
from scipy.optimize import minimize

# 1. Linear PDE
result = solve_linear(PoissonND(), scale=5.0)

# 2. Nonlinear PDE
result = solve_nonlinear(NLPoisson2D(), max_iter=30)

# 3. Inverse: recover source position (x_s,y_s) from sensor data
#    PDE: -Delta u = f,  f = Gaussian at (x_s,y_s),  u=0 on boundary
pde_op = -Op.laplacian(d=2)
basis = SinusoidalBasis.random(2, 800, sigma=5.0, normalize=True)
x_pde = sample_box(3000, 2); x_bc = sample_boundary_box(400, 2)
cache = basis.cache(x_pde)
A_pde = pde_op.apply(basis, x_pde, cache=cache)
A = torch.cat([A_pde, 100*basis.evaluate(x_bc)])
x_sens = torch.tensor([[0.3,0.3],[0.7,0.7],[0.3,0.7],[0.7,0.3]])
def fwd(xs, ys):
    b = torch.exp(-((x_pde[:,0]-xs)**2+(x_pde[:,1]-ys)**2)/0.1).unsqueeze(1)
    b = torch.cat([b, torch.zeros(400,1)])
    beta = solve_lstsq(A, b)
    return (basis.evaluate(x_sens) @ beta).detach().numpy().ravel()
u_obs = fwd(0.4, 0.6) + 0.01*np.random.randn(4)
def loss(p): return np.sum((fwd(p[0], p[1]) - u_obs)**2)
minimize(loss, [0.5, 0.5], method="L-BFGS-B", bounds=[(0.1,0.9)]*2)
\end{lstlisting}

\section{Convergence Profiles}
\label{app:convergence}

\begin{figure}[ht]
\centering
\includegraphics[width=0.48\textwidth]{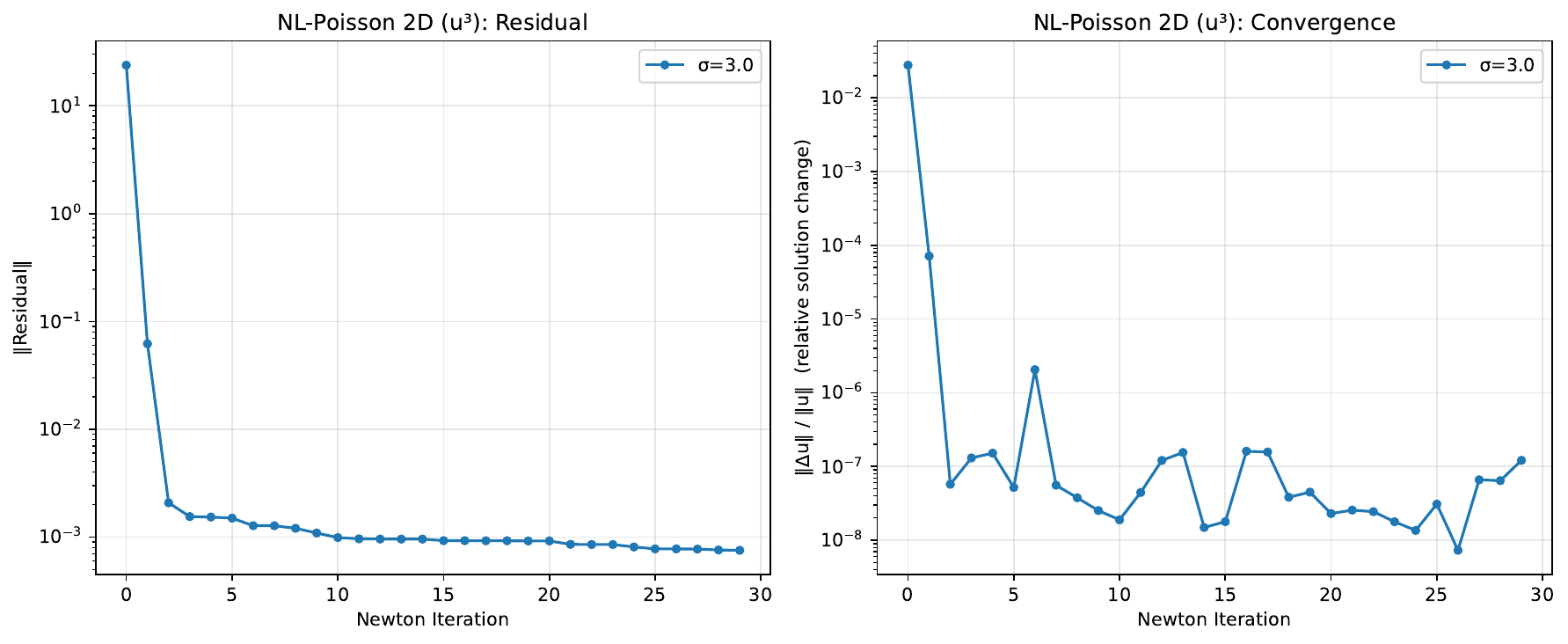}
\includegraphics[width=0.48\textwidth]{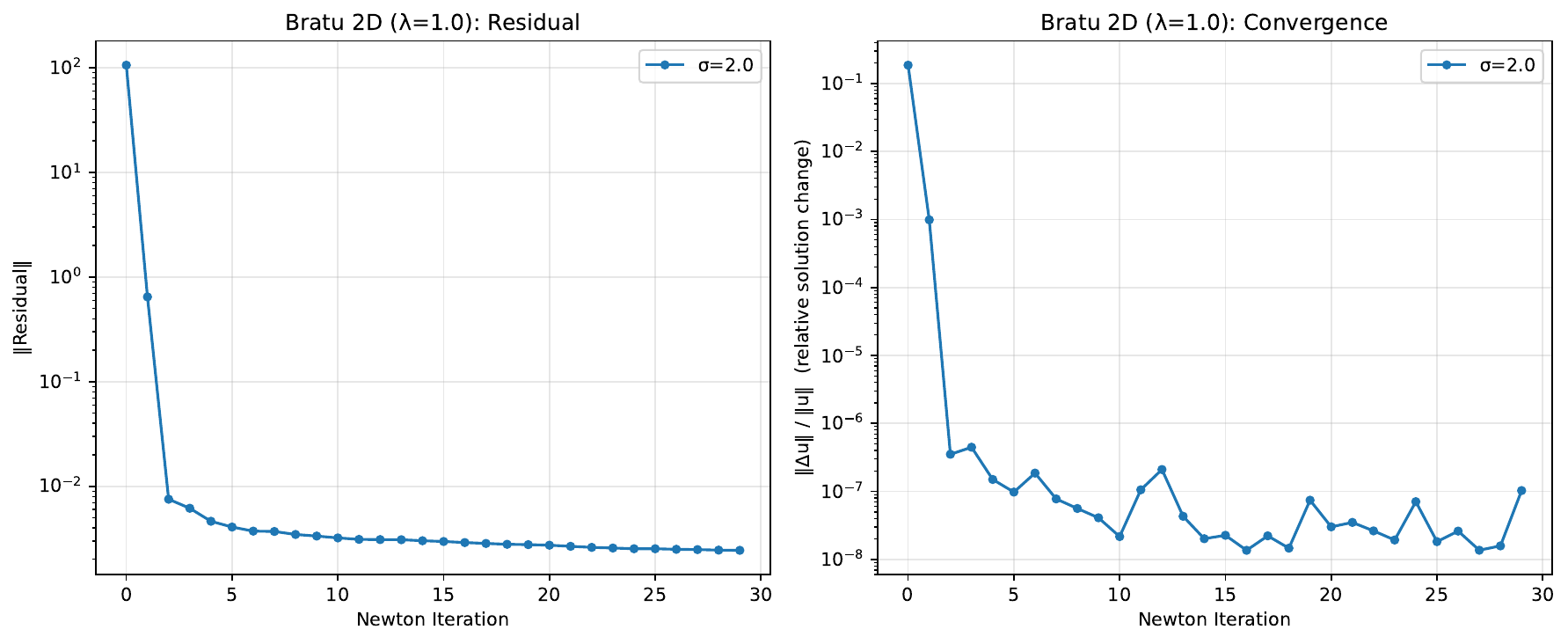}
\caption{Newton convergence: NL-Poisson 2D (left) and Bratu 2D (right).}
\label{fig:conv_nlpoisson}
\end{figure}

\begin{figure}[ht]
\centering
\includegraphics[width=0.48\textwidth]{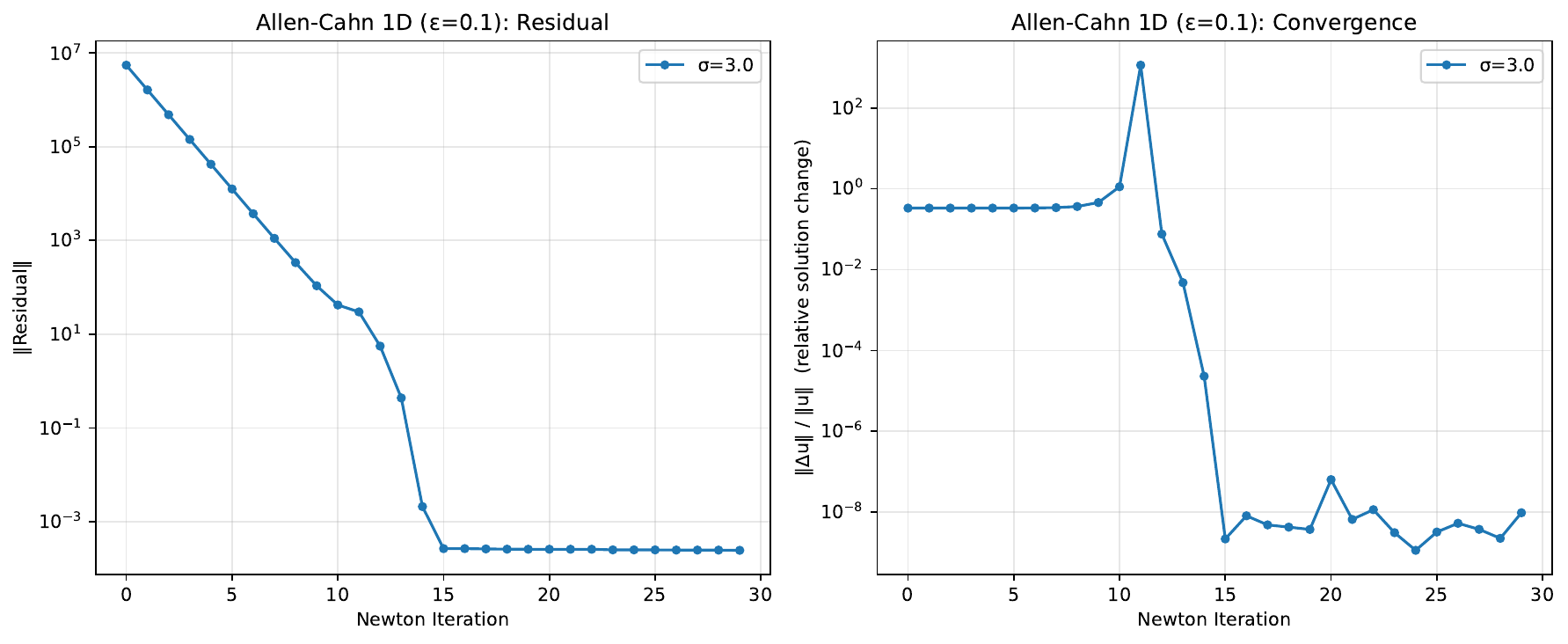}
\caption{Newton convergence: Allen--Cahn 1D.}
\label{fig:conv_allencahn}
\end{figure}

Figures~\ref{fig:conv_nlpoisson}--\ref{fig:conv_allencahn} show Newton convergence profiles (residual norm and relative solution change vs.\ iteration) for for some of the nonlinear problems.
The typical pattern is a rapid initial decrease in residual (2--3 orders of magnitude in the first 2--3 iterations) followed by a slower plateau phase.

\section{Spectral Sensitivity Plots}
\label{app:spectral}

\begin{figure}[ht]
\centering
\includegraphics[width=0.48\textwidth]{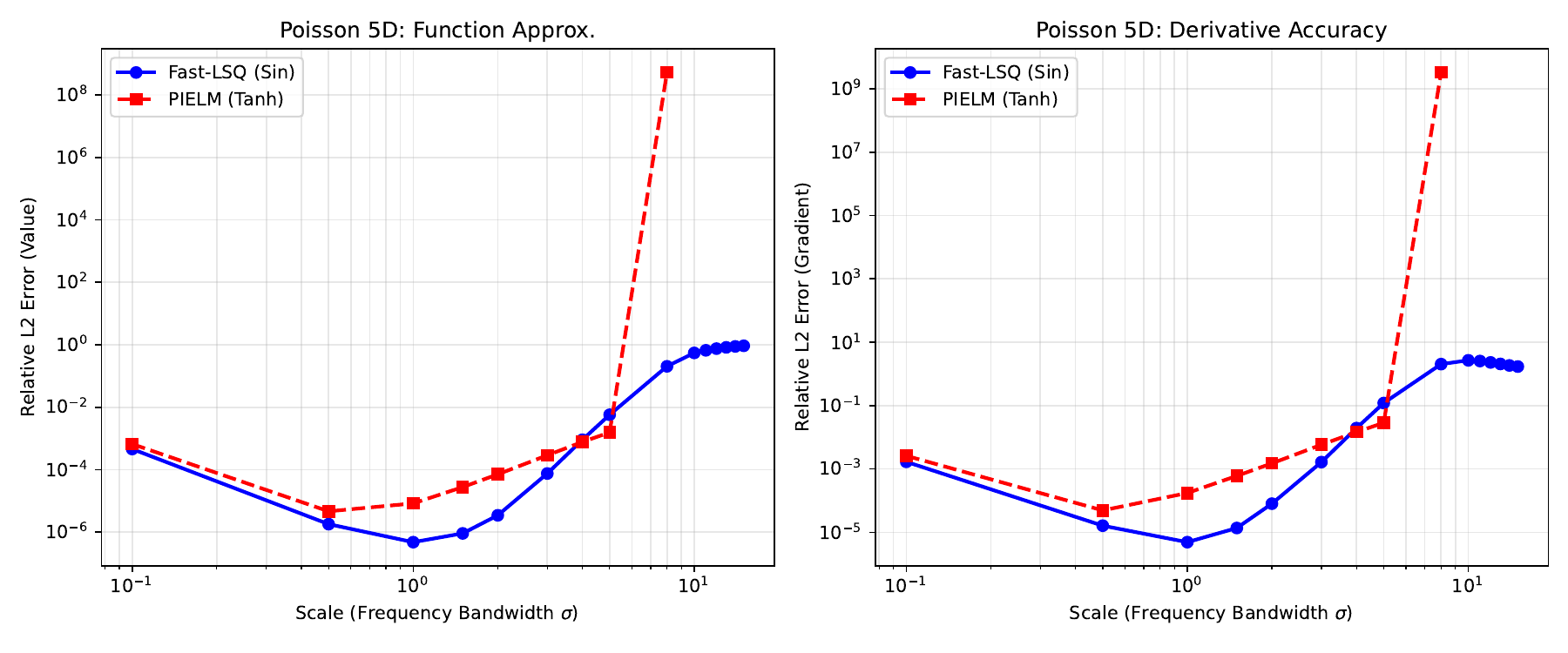}
\includegraphics[width=0.48\textwidth]{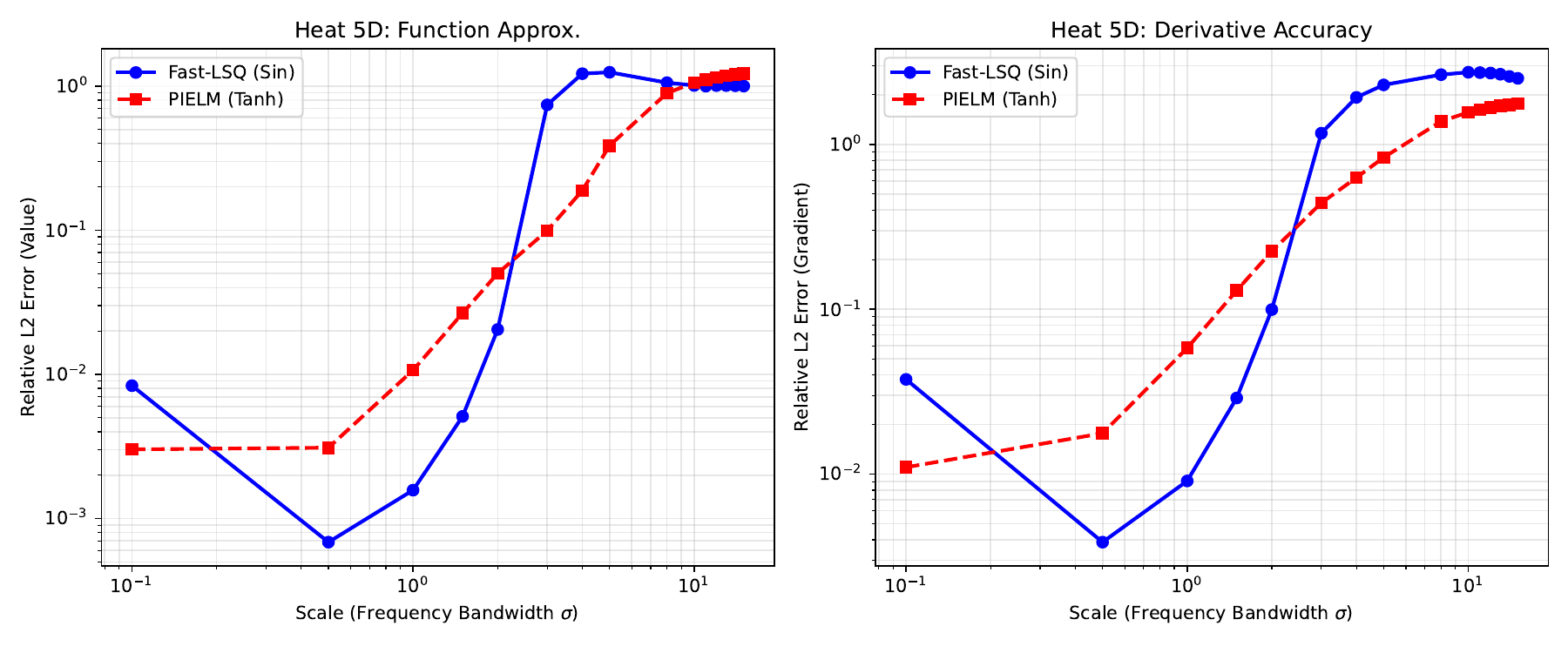}
\caption{Spectral sensitivity: Poisson 5D (left) and Heat 5D (right).}
\label{fig:app_sens1}
\end{figure}

\begin{figure}[ht]
\centering
\includegraphics[width=0.48\textwidth]{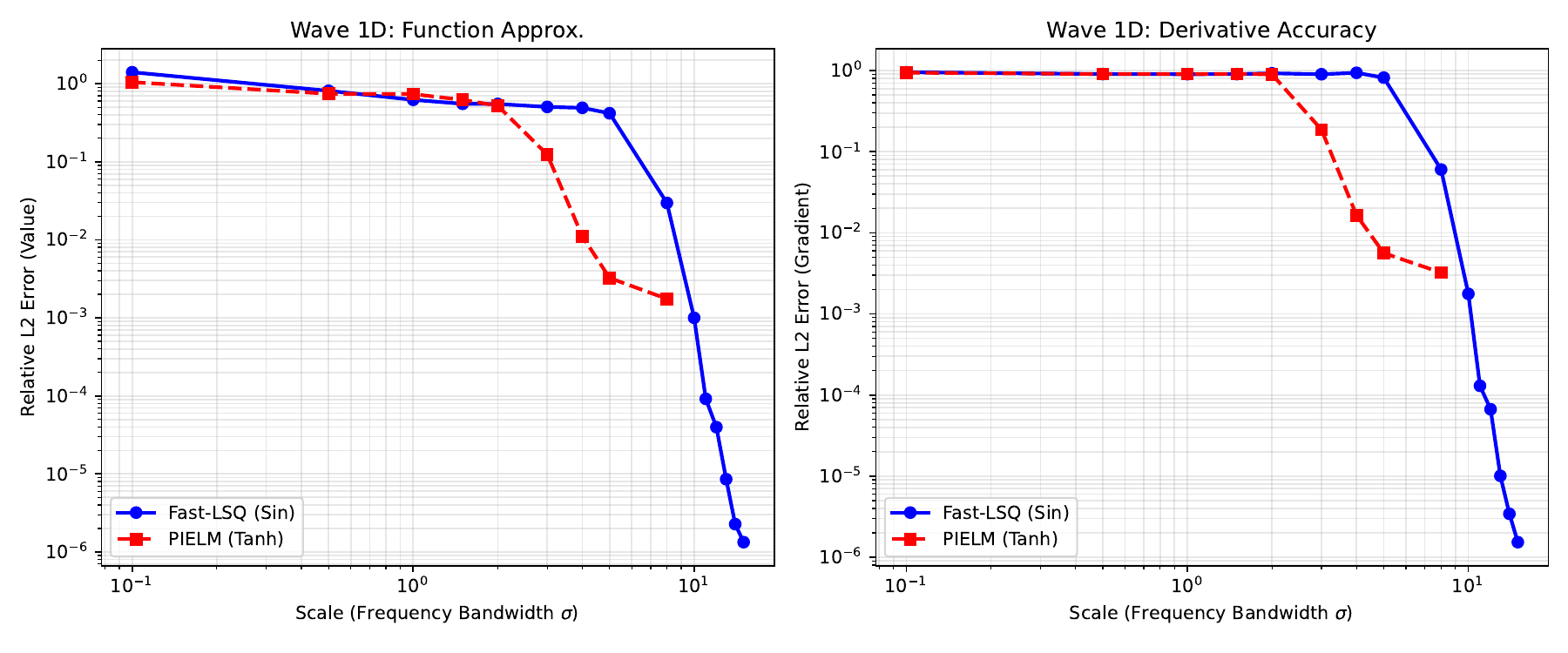}
\includegraphics[width=0.48\textwidth]{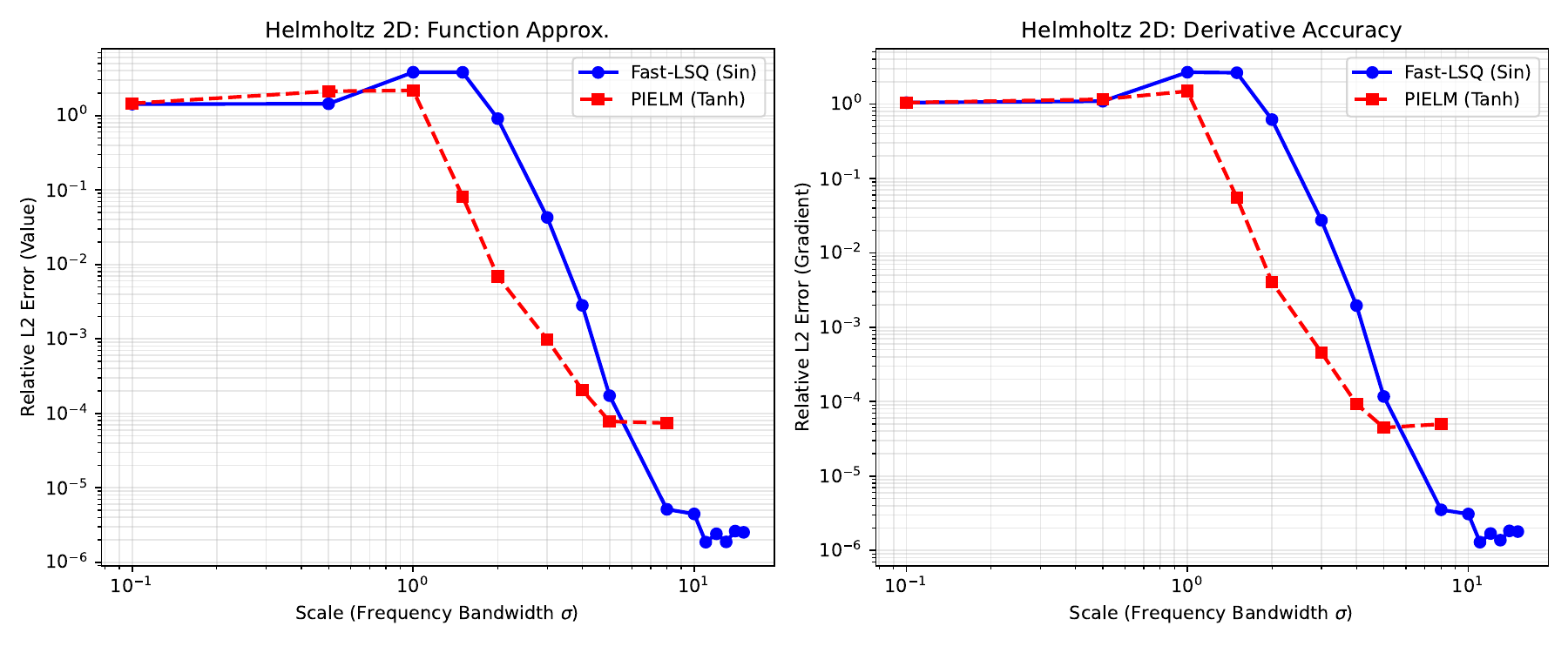}
\caption{Spectral sensitivity: Wave 1D (left) and Helmholtz 2D (right).}
\label{fig:app_sens2}
\end{figure}

\begin{figure}[ht]
\centering
\includegraphics[width=0.48\textwidth]{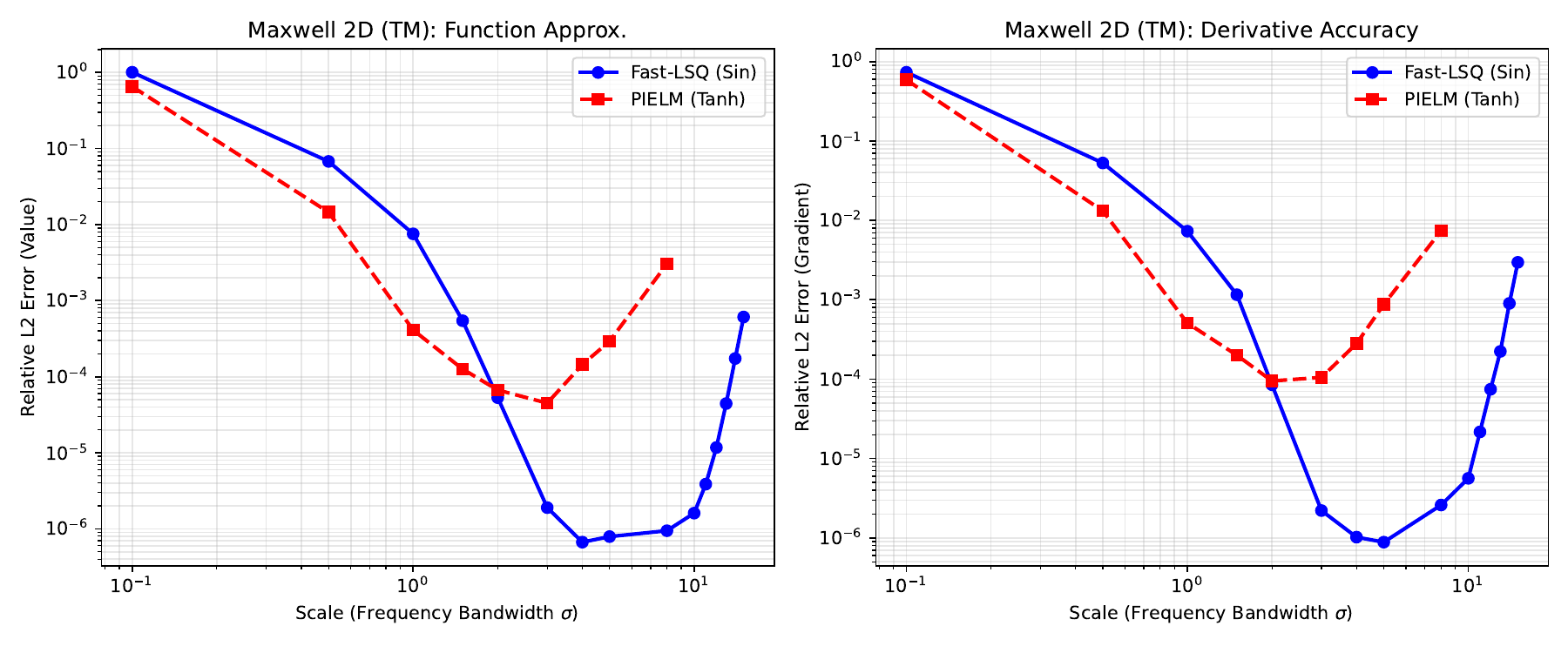}
\includegraphics[width=0.48\textwidth]{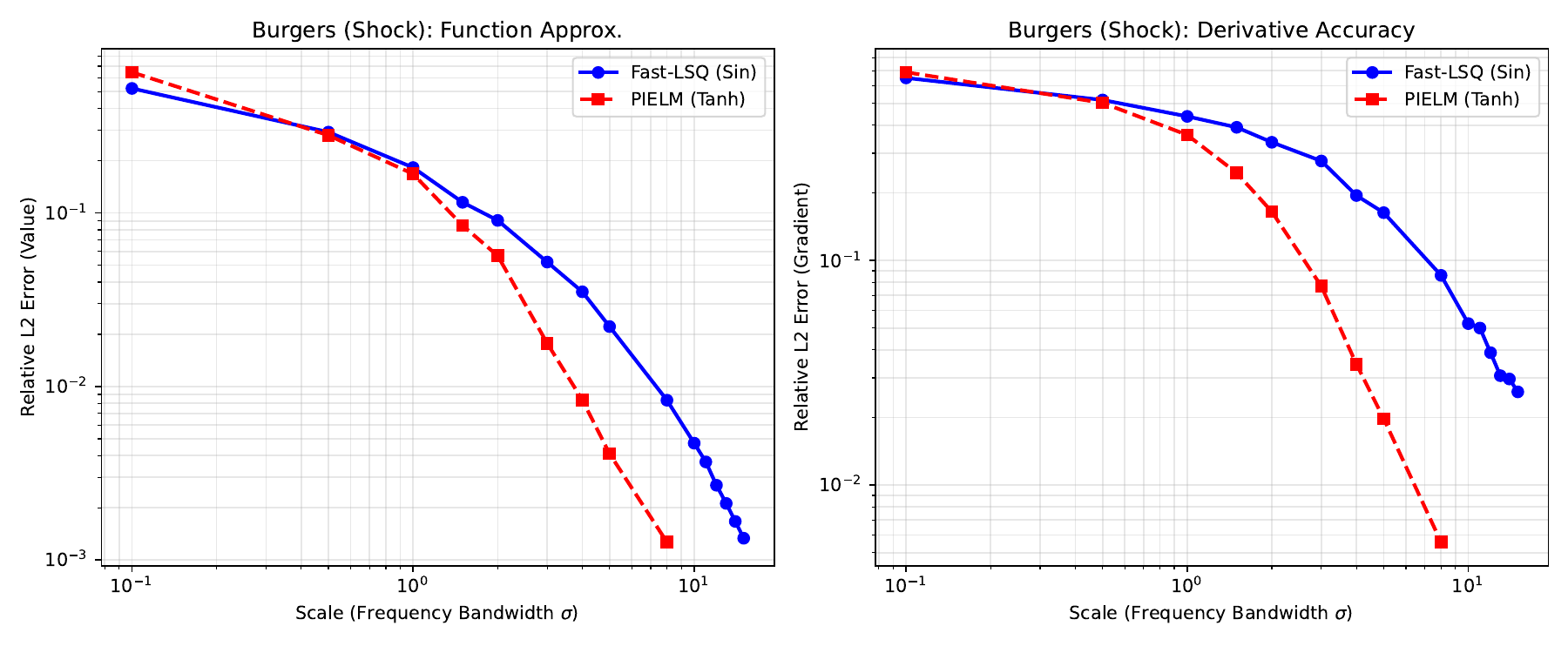}
\caption{Spectral sensitivity: Maxwell 2D (left) and Burgers Shock (right).}
\label{fig:app_sens3}
\end{figure}

\begin{figure}[ht]
\centering
\includegraphics[width=0.48\textwidth]{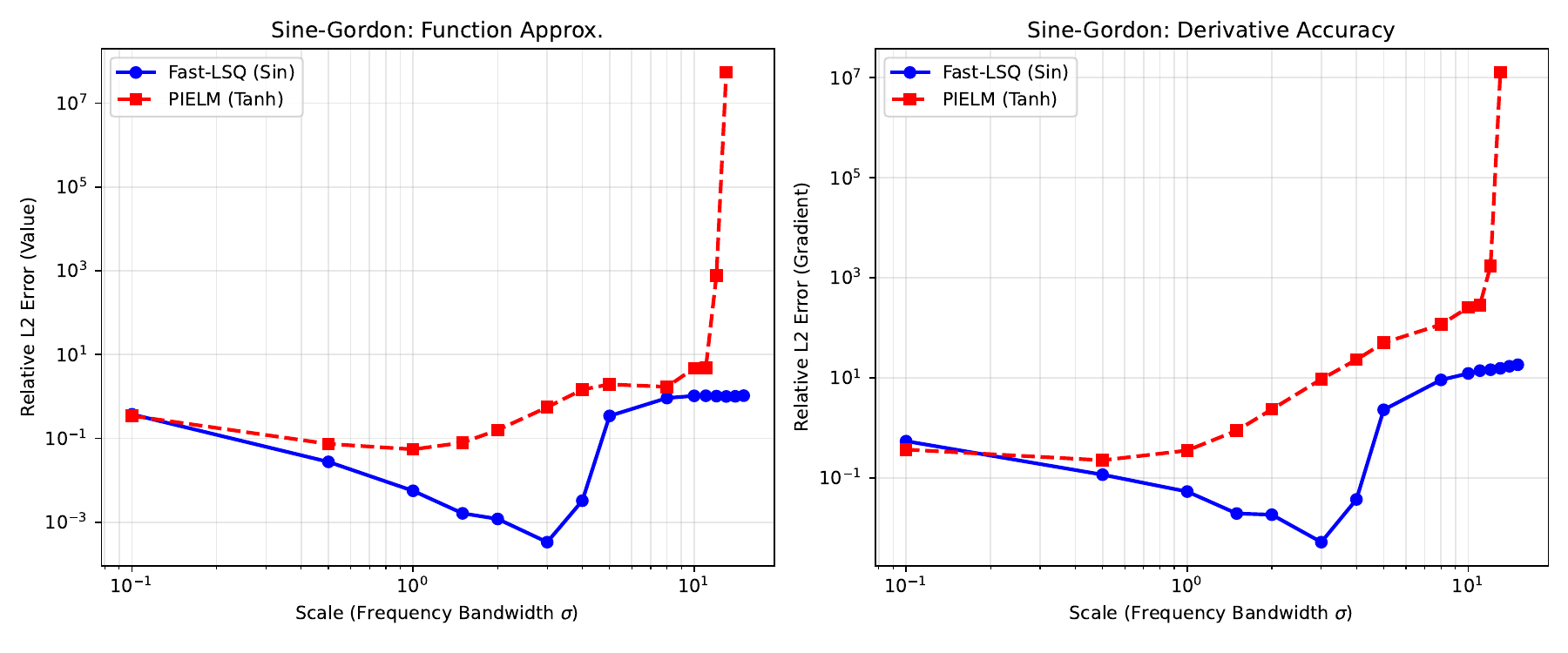}
\includegraphics[width=0.48\textwidth]{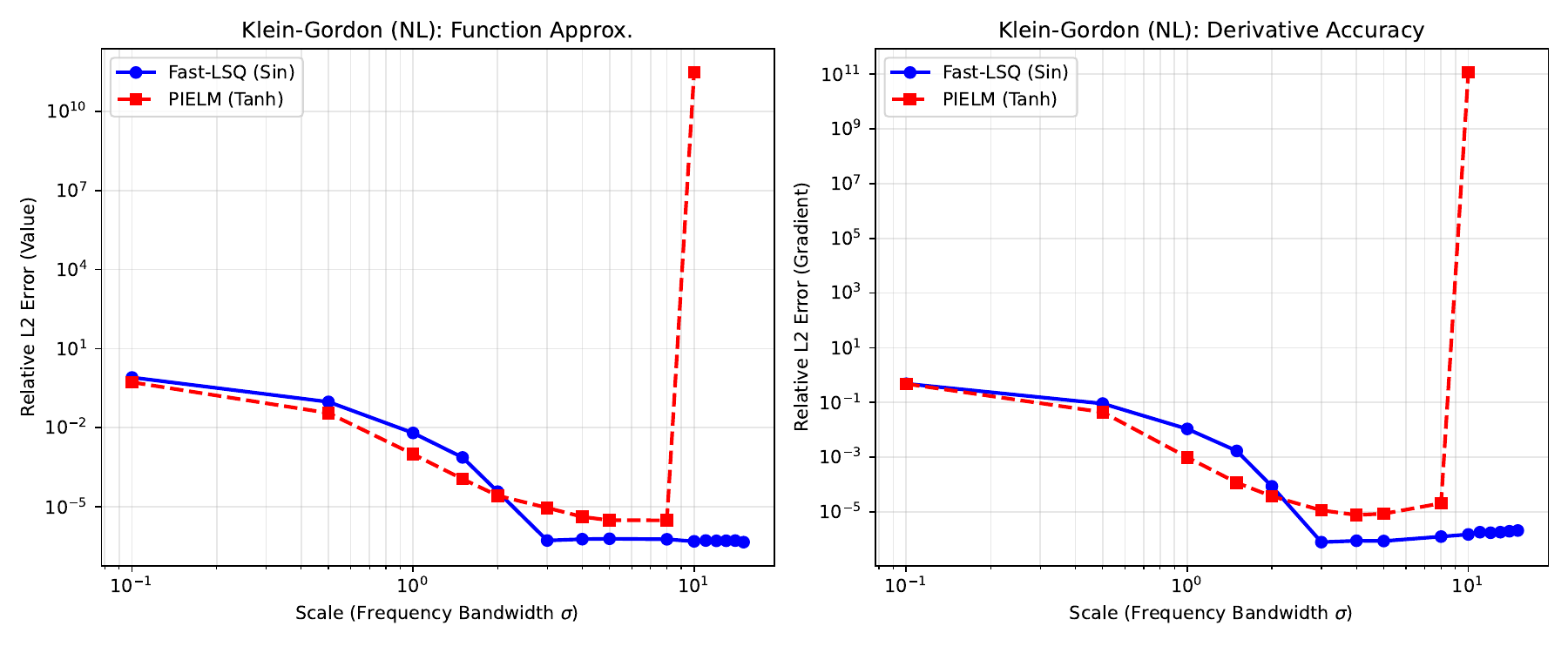}
\caption{Spectral sensitivity: Sine-Gordon (left) and Klein-Gordon (right).}
\label{fig:app_sens4}
\end{figure}

\begin{figure}[ht]
\centering
\includegraphics[width=0.48\textwidth]{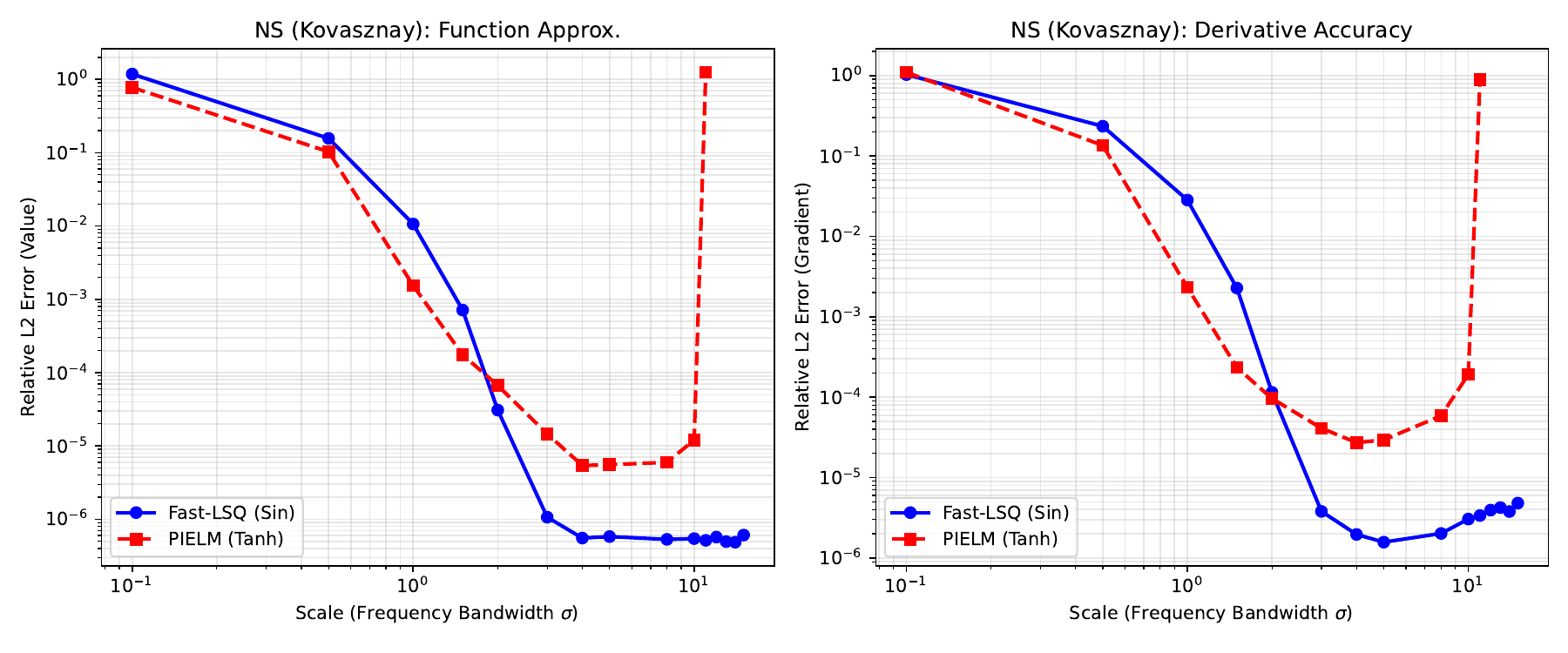}
\includegraphics[width=0.48\textwidth]{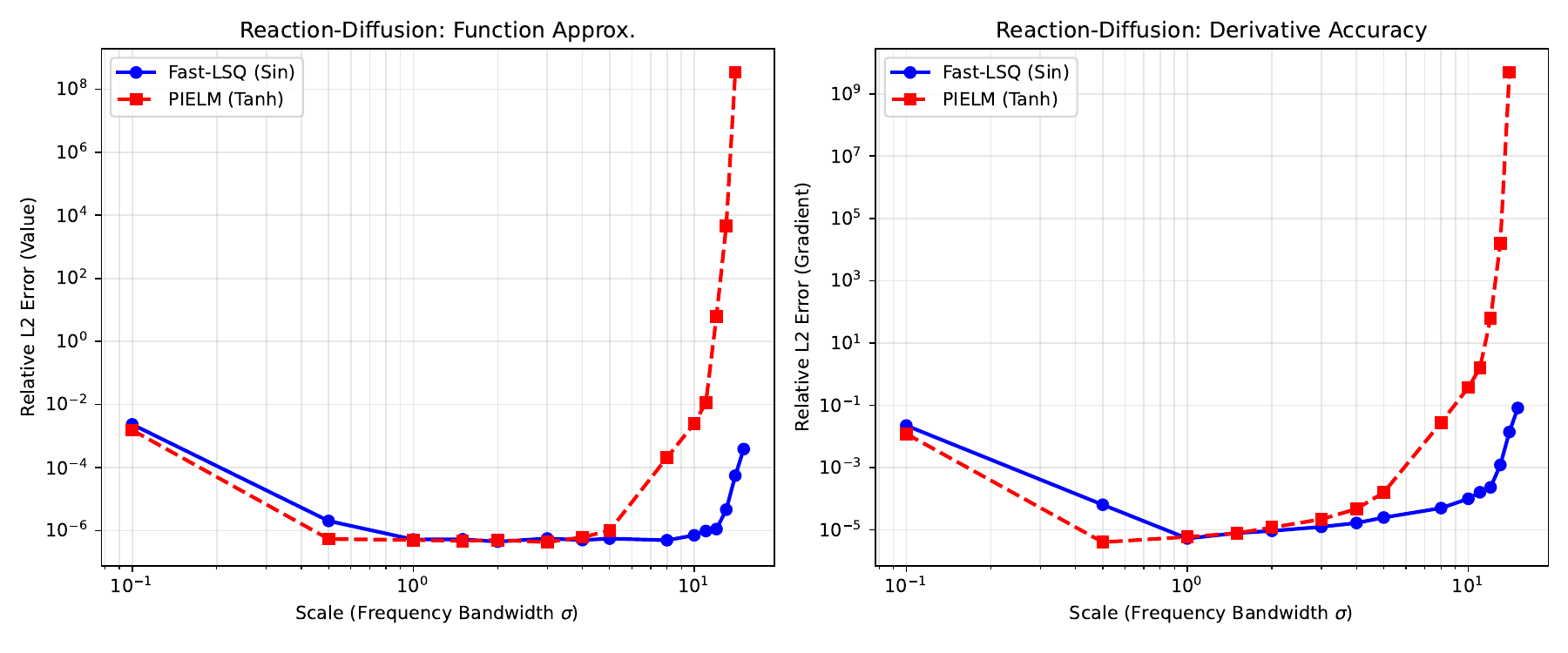}
\caption{Spectral sensitivity: Navier-Stokes Kovasznay (left) and Reaction-Diffusion (right).}
\label{fig:app_sens5}
\end{figure}

\begin{figure}[ht]
\centering
\includegraphics[width=0.48\textwidth]{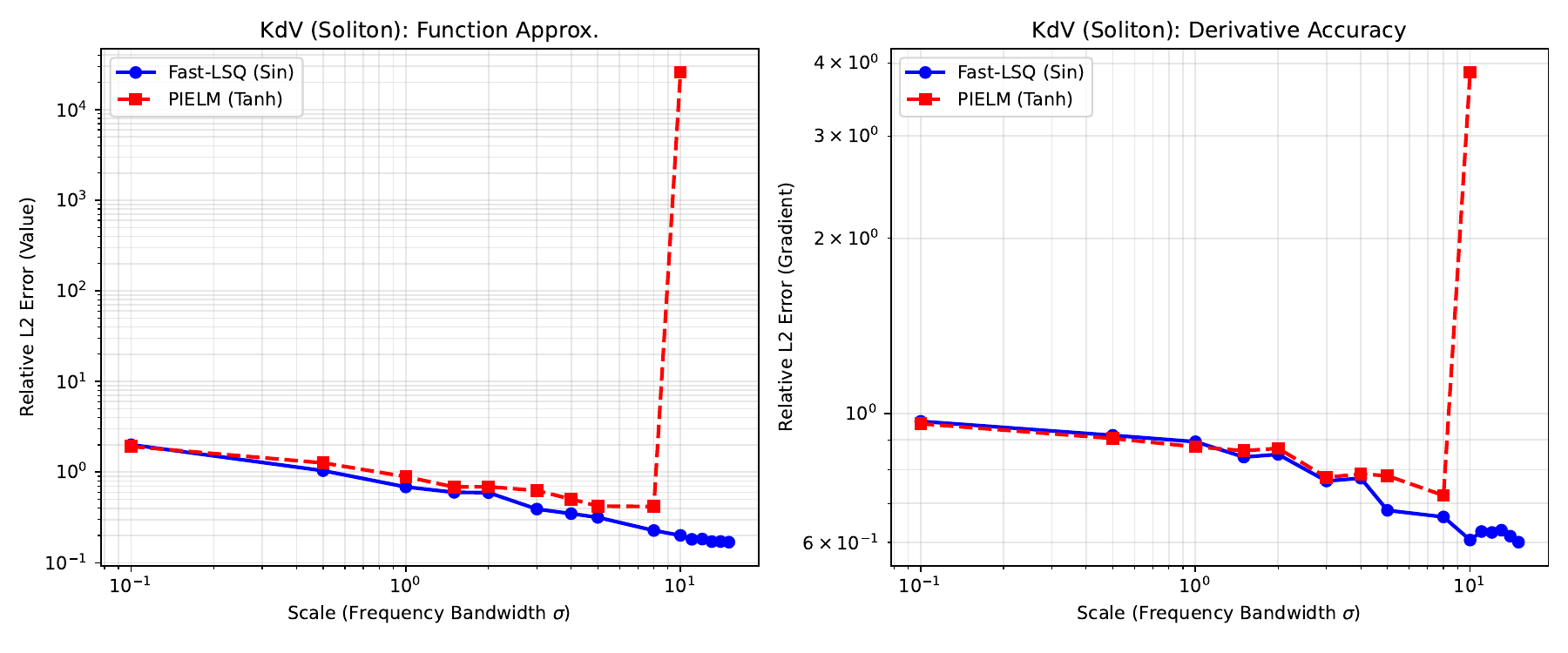}
\includegraphics[width=0.48\textwidth]{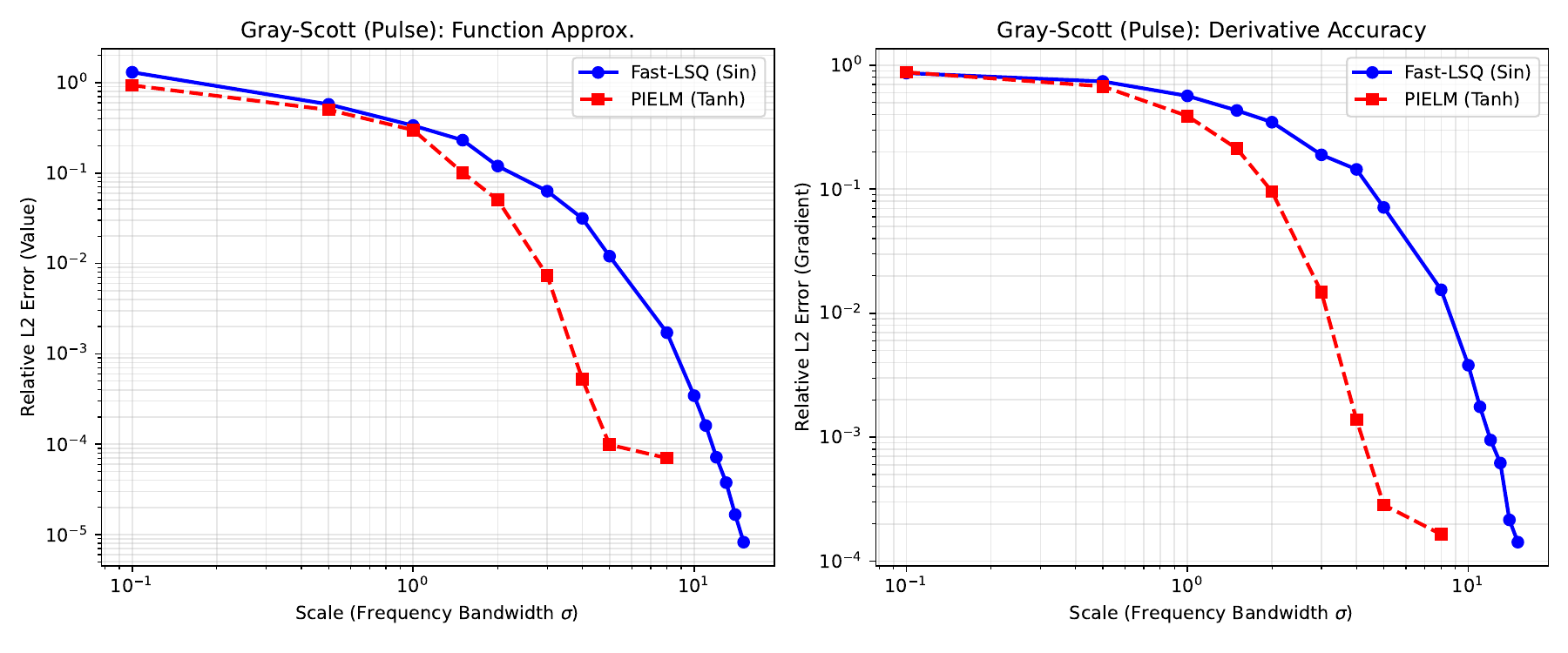}
\caption{Spectral sensitivity: KdV Soliton (left) and Gray-Scott Pulse (right).}
\label{fig:app_sens6}
\end{figure}

Figures~\ref{fig:app_sens1}--\ref{fig:app_sens6} show $L^2$ error vs.\ frequency bandwidth $\sigma$ for all problems.
Each panel shows value error (left) and gradient error (right) for both \fastlsq (sin, blue solid) and PIELM (tanh, red dashed).

\section{Learnable Bandwidth via Reparameterisation}
\label{app:learnable}

This appendix formalises the learnable bandwidth extension outlined in \S\ref{sec:conclusion}.

\vspace{-0.75em}\paragraph{Bandwidth estimation.}
The bandwidth $\sigma$ controls the frequency content of the random Fourier features; selecting it is essential for accuracy.
In practice, we use grid search over $\sigma \in \{0.5, 1, 2, 3, 5, 8, 10, 12, 15\}$ with a few trials per value.
For problems with direction-dependent scales (e.g.\ a wave propagating rapidly in $x$ but diffusing slowly in $y$), isotropic $\sigma$ is suboptimal; an \emph{anisotropic} covariance $\Sigma = LL^\top$ with Cholesky factor $L \in \RR^{d \times d}$ can be learned, giving $d(d{+}1)/2$ parameters that adapt the frequency scale per dimension.
The analytical formulation~\eqref{eq:cyclic} makes both scalar $\sigma$ and full $L$ differentiable through the exact inner solve, so gradient-based optimisation (AdamW or L-BFGS-B) can replace or refine grid search.

\vspace{-0.75em}\subsection{Scalar Bandwidth}

We freeze base weights $\hat{\bW}_j \sim \calN(\mathbf{0}, \mathbf{I}_d)$ and biases $b_j \sim \calU(0, 2\pi)$ at initialisation.
The actual weights are a differentiable function of the learnable parameter $\sigma$:
\begin{equation}
\label{eq:reparam_scalar}
\bW_j(\sigma) = \sigma \hat{\bW}_j, \qquad
\phi_j(\bx;\sigma) = \sin\!\bigl(\sigma \hat{\bW}_j^\top \bx + b_j\bigr).
\end{equation}
Because the inner linear system $\bA(\sigma)\bbeta = \bb$ is solved exactly via least squares, the optimal coefficients are strictly a function of $\sigma$: $\bbeta^*(\sigma) = \bA(\sigma)^\dagger \bb$.
The outer loss becomes
\begin{equation}
\label{eq:outer_loss}
\calL(\sigma) = \bigl\|\bA(\sigma)\bbeta^*(\sigma) - \bb\bigr\|_2^2,
\end{equation}
whose gradient $\nabla_\sigma \calL$ is computable automatically (e.g.\ via \texttt{torch.linalg.lstsq}, which supports backward passes through the pseudo-inverse).
An optimiser such as AdamW~\citep{kingma2014adam} updates $\sigma$ while the coefficients are re-solved exactly at each step, combining the accuracy of spectral feature optimisation with the stability of frozen random features.

\vspace{-0.75em}\subsection{Anisotropic Covariance via Cholesky Decomposition}

For multi-dimensional PDEs with direction-dependent frequency scales (e.g.\ a wave traveling rapidly in $x$ but diffusing slowly in $y$), isotropic sampling $\bW_j \sim \calN(\mathbf{0}, \sigma^2 \mathbf{I}_d)$ is inefficient.
We instead learn a full covariance matrix $\Sigma = LL^\top$, parameterised by a lower-triangular Cholesky factor $L \in \RR^{d \times d}$:
\begin{equation}
\label{eq:reparam_cholesky}
\bW_j(L) = L\hat{\bW}_j, \qquad
\phi_j(\bx;L) = \sin\!\bigl((L\hat{\bW}_j)^\top \bx + b_j\bigr).
\end{equation}
Parameterising via $L$ rather than $\Sigma$ directly guarantees that $\Sigma$ remains symmetric positive semi-definite throughout training.
The number of learnable parameters is $d(d{+}1)/2$, which is modest for typical PDE dimensions ($d \leq 6$).

\vspace{-0.75em}\subsection{Hybrid Training Algorithm}

The resulting procedure is summarised in Algorithm~\ref{alg:learnable}.

\begin{algorithm}[h]
\caption{Learnable-Bandwidth \fastlsq}
\label{alg:learnable}
\begin{algorithmic}[1]
\Require PDE operator $\calL$, BC operator $\calB$, source $f$, BC data $g$
\Require Number of features $N$, initial Cholesky factor $L_0$ (e.g.\ $\mathbf{I}_d$)
\State Freeze base weights $\hat{\bW}_j \sim \calN(\mathbf{0}, \mathbf{I}_d)$, $b_j \sim \calU(0, 2\pi)$, $j = 1,\ldots,N$
\While{not converged (AdamW outer steps)}
    \State $\bW_j \gets L \hat{\bW}_j$ \Comment{Reparameterise}
    \State Assemble $\bA(L)$ analytically via~\eqref{eq:cyclic}
    \State $\bbeta^* \gets \bA(L)^\dagger \bb$ \Comment{Inner exact solve}
    \State $\calL_{\text{outer}} \gets \|\bA(L)\bbeta^* - \bb\|_2^2$ \Comment{Outer loss}
    \State $L \gets L - \eta\,\nabla_L \calL_{\text{outer}}$ \Comment{AdamW update}
\EndWhile
\State \Return $u_N(\bx) = \frac{1}{\sqrt{N}}\sum_j \beta_j^* \sin((L\hat{\bW}_j)^\top \bx + b_j)$
\end{algorithmic}
\end{algorithm}

This hybrid approach retains the extreme accuracy of the one-shot exact solve for $\bbeta$ while using gradient-based optimisation only for the low-dimensional bandwidth parameters, avoiding the chaotic, high-variance gradients that arise when training the individual weight entries $\bW_j$ of a standard PINN.

\section{Trigonometric Symmetries Exploited by \fastlsq}
\label{app:symmetries}

The effectiveness of sinusoidal random features for PDE solving rests on a rich set of trigonometric identities.
This appendix provides a systematic taxonomy, organised by the role each identity plays in the \fastlsq framework.

\vspace{-0.75em}\subsection{Calculus Symmetries (Foundation)}

These describe how sinusoids behave under differentiation and form the core mechanism that eliminates automatic differentiation.

\vspace{-0.75em}\paragraph{Cyclic derivative (modulo 4).}
The $n$-th derivative cycles through a fixed four-element sequence:
\begin{equation}
\frac{d^n}{dz^n}\sin(z) = \Phi_{n \bmod 4}(z), \qquad
\Phi_0 = \sin,\; \Phi_1 = \cos,\; \Phi_2 = -\sin,\; \Phi_3 = -\cos.
\end{equation}

\vspace{-0.75em}\paragraph{Chain rule pull-out.}
For a linear combination $z = \bW^\top \bx + b$ with multi-index $\alpha$:
\begin{equation}
D^\alpha \sin(\bW^\top \bx + b) = \Bigl(\prod_k W_k^{\alpha_k}\Bigr)\, \Phi_{|\alpha| \bmod 4}(\bW^\top \bx + b).
\end{equation}
This is precisely equation~\eqref{eq:cyclic} of the main text.

\vspace{-0.75em}\paragraph{Eigenfunction property.}
Sinusoids are eigenfunctions of the Laplacian:
\begin{equation}
\Delta \sin(\bW^\top \bx + b) = -\|\bW\|^2 \sin(\bW^\top \bx + b),
\end{equation}
and more generally of any constant-coefficient linear differential operator.
This makes operator evaluation a scalar multiplication, which is why Corollary~\ref{cor:operators} holds.

\vspace{-0.75em}\subsection{Algebraic Phase and Shift Symmetries}

These identities enable algebraic manipulation of the inner weights and biases.

\vspace{-0.75em}\paragraph{Angle addition (phase decoupling).}
$\sin(x + y) = \sin(x)\cos(y) + \cos(x)\sin(y)$.
This allows splitting a shifted feature $\sin(\bW^\top\bx + b)$ into bias-free sine and cosine components weighted by $\cos(b)$ and $\sin(b)$, which can simplify analysis of the feature distribution.

\vspace{-0.75em}\paragraph{Parity (even/odd symmetry).}
$\sin(-x) = -\sin(x)$ and $\cos(-x) = \cos(x)$.
Useful for hard-coding symmetrical physical boundary conditions: if $u(-x) = u(x)$, one can restrict to even-parity features (cosines only).

\vspace{-0.75em}\paragraph{Quarter-wave and half-wave shifts.}
$\sin(x + \tfrac{\pi}{2}) = \cos(x)$, $\cos(x - \tfrac{\pi}{2}) = \sin(x)$, and $\sin(x + \pi) = -\sin(x)$.
These relate bias shifts to amplitude sign flips and sine-cosine swaps, explaining why the bias distribution $b_j \sim \calU(0, 2\pi)$ provides sufficient coverage of both function families.

\vspace{-0.75em}\subsection{Product and Power Identities (Nonlinearity Linearisation)}

When dealing with nonlinear PDEs whose residuals contain products or powers of $u$ (e.g.\ the $u^3$ term in NL-Poisson), these identities show that products of sinusoids remain within the sinusoidal family.

\vspace{-0.75em}\paragraph{Product-to-sum.}
$\sin(x)\sin(y) = \tfrac{1}{2}[\cos(x - y) - \cos(x + y)]$. Multiplying two random features produces two new features at sum and difference frequencies.

\vspace{-0.75em}\paragraph{Power reduction.}
$\sin^2(x) = \tfrac{1 - \cos(2x)}{2}$ and $\cos^2(x) = \tfrac{1 + \cos(2x)}{2}$.
Squaring doubles the frequency.

\vspace{-0.75em}\paragraph{Cubic reduction.}
$\sin^3(x) = \tfrac{3\sin(x) - \sin(3x)}{4}$.
This identity is directly relevant to the $u^3$ residual term in the NL-Helmholtz benchmark (\S\ref{sec:newton_results}): evaluating $u^3$ on a sinusoidal expansion analytically generates features at the base and triple frequencies, remaining within the sinusoidal family without point-wise grid sampling.

\vspace{-0.75em}\subsection{Integral and Orthogonality Symmetries}

These properties are relevant for global (integral-form) loss evaluation and theoretical analysis.

\vspace{-0.75em}\paragraph{Orthogonality.}
$\int_{-\pi}^{\pi} \sin(nx)\sin(mx)\,dx = 0$ for $n \neq m$, and $\int_{-\pi}^{\pi} \sin^2(nx)\,dx = \pi$.
This underpins the spectral accuracy of Fourier-type expansions: distinct frequency components are orthogonal, so the least-squares solve decomposes cleanly in the frequency domain when collocation points are dense enough to approximate the integral.

\vspace{-0.75em}\subsection{Connection to Euler's Formula}

All identities above unify under the complex exponential representation:
\begin{equation}
e^{ix} = \cos(x) + i\sin(x), \qquad
\sin(x) = \frac{e^{ix} - e^{-ix}}{2i}, \qquad
\cos(x) = \frac{e^{ix} + e^{-ix}}{2}.
\end{equation}
Every trigonometric identity reduces to elementary exponent rules ($e^a e^b = e^{a+b}$).
This perspective also connects \fastlsq to the Fourier transform of the Gaussian RBF kernel (\S\ref{sec:rff}): the random feature approximation $\frac{1}{N}\sum_j \phi_j(\bx)\phi_j(\bx')$ converges to $\exp(-\frac{\sigma^2}{2}\|\bx - \bx'\|^2)$ precisely because the Fourier transform of the Gaussian is Gaussian, and each sinusoidal feature corresponds to a single frequency sample from this transform.

\vspace{-0.75em}\subsection{Theoretical Guarantees}

\vspace{-0.75em}\paragraph{Universal approximation.}
The Gaussian RBF kernel spans a Reproducing Kernel Hilbert Space (RKHS) that is dense in the space of continuous functions on compact sets.
The $1/\sqrt{N}$-scaled random features converge to this kernel as $N \to \infty$~\citep{rahimi2007random}, guaranteeing that given enough features, \fastlsq can approximate any sufficiently smooth PDE solution to arbitrary precision.

\vspace{-0.75em}\paragraph{Convergence rate.}
Standard random feature approximation error scales as $\calO(1/\sqrt{N})$~\citep{rahimi2007random}, with problem-dependent constants related to the smoothness of the target function in the RKHS norm.

\vspace{-0.75em}\paragraph{Conditioning limits.}
As $N$ grows very large, the feature matrix $\bA$ becomes increasingly ill-conditioned because the columns become near-linearly-dependent.
For high-order PDEs, the monomial prefactor $\prod W_k^{\alpha_k}$ in~\eqref{eq:cyclic} amplifies the condition number further, explaining the accuracy plateau observed in some experiments (\S\ref{sec:spectral}).
The Tikhonov regularisation parameter $\mu$ mitigates this at the cost of a small bias.

\section{Extension Results}
\label{app:extensions}

This appendix provides additional details for the extensions in \S\ref{sec:conclusion} and \S\ref{sec:applications}.
All experiments run on a single CPU core.

\vspace{-0.75em}\subsection{PDE Discovery Details}
\label{app:pde_discovery}

The PDE discovery experiment of \S\ref{sec:applications} uses $M{=}2000$ data points, $N{=}1500$ features ($\sigma{=}3$), Tikhonov $\mu{=}10^{-8}$, noise level $\sigma_\varepsilon{=}0.01$, and a LASSO sweep over $\lambda \in [10^{-5}, 10^1]$ (100 values).
After LASSO selects the equation structure $\{u, u_x\}$, an OLS refit on the active terms removes shrinkage bias.
The recovered coefficients are $c_1 = -4.26$ and $c_2 = -1.12$ (true: $-4.25$ and $-1.00$); the stiffness is recovered to ${<}0.2\%$ relative error while the damping term has $\sim\!12\%$ error, reflecting its smaller contribution to $u_{xx}$ ($|c_2|/|c_1| \approx 0.24$).
Finite-difference $u_{xx}$ at this noise level has RMSE~$\approx 2461$, rendering standard SINDy inoperable; the analytical derivatives achieve RMSE~$\approx 0.41$, a $\sim\!6000\times$ improvement.

\vspace{-0.75em}\subsection{Inverse Heat Source Details}
\label{app:inverse_heat}

The inverse heat-source problem of \S\ref{sec:applications} uses the 2-D heat equation
$u_t - 0.05\,\Delta u = \sum_{k=1}^{4} f_k(\bx;\boldsymbol{\theta}_k)$
on $[0,1]^2 \times [0, 1.5]$ with zero Dirichlet BCs and zero initial condition.
Each source $f_k$ is an anisotropic Gaussian with precision matrix
$P_k = L_k^\top L_k$, $L_k = \bigl[\begin{smallmatrix}a_k & 0 \\ b_k & c_k\end{smallmatrix}\bigr]$,
so $\boldsymbol{\theta}_k = (x_s, y_s, I, a, b, c)$ gives 24 parameters in total.
The operator matrix (1200 features, 6000 interior + 2000 BC + 1000 IC points)
is assembled once and Cholesky-prefactored; each of the 150 L-BFGS-B iterations
requires only two triangular back-substitutions.
Table~\ref{tab:inverse_heat} gives the recovery results for source positions and intensities
(the primary quantities of interest; shape parameters are weakly constrained and omitted).
Only 4 irregularly placed sensors are used, each recording a temperature time-series
over 60 snapshots (240 observations total) as the space-time field $u(x,y,t)$ evolves.
Three sources are localised to within $0.02$ per coordinate; the weakest source (Source~4)
has a larger position error ($0.12$), reflecting the genuine information limit of 4 sensors
for a 24-parameter problem.

\begin{table}[ht]
\caption{Inverse recovery of 4 heat-source positions and intensities
(1200 features, 4 sensors $\times$ 60 time snapshots = 240 observations, 300 L-BFGS-B iterations).
Shape parameters ($a$, $b$, $c$) are only weakly constrained by 4 sensors; position recovery
is the primary target.}
\label{tab:inverse_heat}
\vskip 0.05in
\centering\small
\setlength{\tabcolsep}{5pt}
\begin{tabular}{l rrr rrr rrr rrr}
\toprule
 & \multicolumn{3}{c}{Source 1 (BL)} & \multicolumn{3}{c}{Source 2 (BR)}
 & \multicolumn{3}{c}{Source 3 (TR)} & \multicolumn{3}{c}{Source 4 (TL)} \\
\cmidrule(lr){2-4} \cmidrule(lr){5-7} \cmidrule(lr){8-10} \cmidrule(lr){11-13}
Param & True & Rec. & $|\Delta|$ & True & Rec. & $|\Delta|$ & True & Rec. & $|\Delta|$ & True & Rec. & $|\Delta|$ \\
\midrule
$x_s$ & 0.22 & 0.221 & 9e-4 & 0.78 & 0.847 & 6.7e-2 & 0.75 & 0.768 & 1.8e-2 & 0.25 & 0.131 & 1.2e-1 \\
$y_s$ & 0.22 & 0.211 & 8.9e-3 & 0.25 & 0.291 & 4.1e-2 & 0.75 & 0.759 & 9.3e-3 & 0.78 & 0.722 & 5.8e-2 \\
$I$   & 4.00 & 4.26  & 2.6e-1 & 2.50 & 3.37  & 8.7e-1 & 3.50 & 3.98  & 4.8e-1 & 2.00 & 2.55  & 5.5e-1 \\
\bottomrule
\end{tabular}
\vskip -0.1in
\end{table}

\vspace{-0.75em}\subsection{Inverse Magnetostatics Details}
\label{app:inverse_magneto}

The coil localisation setup of \S\ref{sec:applications} uses $N{=}1500$ features, $M_{\rm int}{=}6000$, $M_{\rm bc}{=}1500$, $\lambda{=}200$.
The iron geometry has four hyperbolic pole faces ($|x^2 - y^2| = R_{\rm ap}^2$, $R_{\rm ap}{=}0.35$) with $\mu_r(x,y) \in [1, 50]$.
Eight sensors on a circle of radius 0.25 observe $B_x$, $B_y$ with 1\% noise.
We minimise $\calL(x_c, y_c) = \frac{1}{8}\sum_{s=1}^{8} [(B_x - \hat{B}_x^s)^2 + (B_y - \hat{B}_y^s)^2]$ via 40 Adam steps ($\eta{=}0.04$).

\vspace{-0.75em}\subsection{Comprehensive Baseline Comparison: All 10 PDEs}
\label{app:rbf_fd}

To complement the PIELM and PINNacle comparisons in \S\ref{sec:experiments},
we benchmark \fastlsq against two conventional solver families on \emph{all
10 PDEs from the paper}.  This stress-tests applicability as well as accuracy.

Table~\ref{tab:all_equations} lists all equations we run, their PDE form, domain, and the scikit-fem P2 FEM benchmark where applicable.
scikit-fem applies only to 2-D elliptic spatial problems (Helmholtz, NL-Poisson, Bratu, NL-Helmholtz); Poisson~5D and Heat~5D are $d{\ge}5$ (FEM mesh scales as $h^{-d}$, intractable); Wave~1D and Maxwell~2D~TM are hyperbolic space-time; Burgers~1D and Allen--Cahn~1D are 1-D steady BVPs solved by scipy.solve\_bvp.

\begin{table}[t]
\caption{All 10 PDEs benchmarked: equation, domain, \fastlsq $L^2$ error, and conventional baseline. scikit-fem P2 FEM (refine 4--6, ${\approx}4$k--17k DoF) for 2-D elliptic; scipy.solve\_bvp for 1-D steady BVPs. ``N/A'' = conventional solver not applicable.}
\label{tab:all_equations}
\vskip 0.02in
\centering\scriptsize
\begin{tabular}{ll cc c}
\toprule
Problem & PDE & Domain & \fastlsq $L^2$ & scikit-fem / solve\_bvp \\
\midrule
\multicolumn{5}{l}{\textit{Linear}} \\
Poisson 5D   & $-\Delta u = f$ & $[0,1]^5$ & 4.8e-7 & N/A ($d{=}5$) \\
Heat 5D      & $u_t - \kappa\Delta u = f$ & $B^5{\times}[0,1]$ & 6.9e-4 & N/A ($d{=}6$) \\
Wave 1D      & $u_{tt} - c^2 u_{xx} = 0$ & $[0,1]^2$ (space-time) & 1.3e-6 & N/A (hyperbolic) \\
Helmholtz 2D & $\Delta u + k^2 u = f$ & $[0,1]^2$ & 1.9e-6 & 4.0e-5 \\
Maxwell 2D   & $u_{tt} - c^2(u_{xx}{+}u_{yy}) = 0$ & $[0,1]^3$ (space-time) & 6.7e-7 & N/A (hyperbolic) \\
\midrule
\multicolumn{5}{l}{\textit{Nonlinear}} \\
NL-Poisson 2D & $-\Delta u + u^3 = f$ & $[0,1]^2$ & 6.1e-8 & 2.6e-6 \\
Bratu 2D     & $-\Delta u - \lambda e^u = f$ & $[0,1]^2$ & 1.5e-7 & 2.6e-6 \\
Burgers 1D   & $u u' - \nu u'' = f$ & $[0,1]$ & 3.9e-9 & 1.8e-10 (solve\_bvp) \\
NL-Helm. 2D  & $\Delta u + k^2 u + \alpha u^3 = f$ & $[0,1]^2$ & 2.4e-9 & 2.4e-6 \\
Allen--Cahn 1D & $\varepsilon u'' + u - u^3 = f$ & $[0,1]$ & 6.0e-8 & 1.2e-10 (solve\_bvp) \\
\bottomrule
\end{tabular}
\vskip -0.05in
\end{table}

\vspace{-0.75em}\paragraph{RBF Kansa collocation.}
We implement the Kansa method~\citep{kansa1990multiquadrics} using Hardy
multiquadric (MQ) radial basis functions $\phi_j(\bx){=}\sqrt{c^2{+}|\bx{-}\mathbf{c}_j|^2}$,
$c{=}0.5$, placed at $N{\le}1500$ scattered centres.
The Laplacian $\Delta\phi_j = (d\,c^2 + (d{-}1)\,r^2)\phi_j^{-3}$ is assembled
analytically in any dimension, giving a one-shot overdetermined system.
Newton--Raphson with backtracking is used for nonlinear problems.
RBF is dimension-agnostic in principle but exhibits poor conditioning for
high-frequency ($k{=}10$) or stiff 1-D problems with thin layers.

\vspace{-0.75em}\paragraph{scikit-fem P2 FEM (2-D elliptic problems).}
We use the \texttt{scikit-fem}~\citep{scikit-fem} package~(pip-installable, pure Python)
with P2 Lagrange triangular elements on a uniform mesh refined to
${\approx}4{,}000$--$17{,}000$ DoF.
Newton--Raphson with full Galerkin Jacobian assembly handles nonlinear problems.
FEM is restricted to \emph{2-D elliptic spatial} problems; high-dimensional
($d{\ge}5$) and hyperbolic space-time problems are outside its scope.

\vspace{-0.75em}\paragraph{scipy.integrate.solve\_bvp (1-D steady BVPs).}
For 1-D steady BVPs (Burgers, Allen-Cahn) we use scipy's built-in
adaptive collocation solver (equivalent to MATLAB's \texttt{bvp4c}).
This method is not applicable to 2-D or higher-dimensional problems.

\vspace{-0.75em}\paragraph{Summary.}
Across all 10 PDEs, \fastlsq is the \emph{only} method that produces valid results for every equation: the high-dimensional problems (Poisson~5D, Heat~5D) and hyperbolic space-time problems (Wave~1D, Maxwell~2D~TM) lie entirely outside the scope of conventional grid-based solvers, while the 2-D elliptic and 1-D BVP problems are solved by the respective specialised methods but at lower accuracy than \fastlsq (see Table~\ref{tab:all_equations} and the Conv.\ columns in Tables~\ref{tab:linear_results} and~\ref{tab:newton_results}).
On the 2-D nonlinear problems, scikit-fem P2 FEM achieves $\approx2.6\times10^{-6}$ at 4000~DoF while \fastlsq reaches $10^{-7}$--$10^{-9}$ at 1500 features.
RBF Kansa diverges on the stiff Burgers problem and loses two orders of magnitude relative to \fastlsq on all 2-D nonlinear problems, confirming that the sinusoidal basis is more numerically robust than the MQ-RBF basis at equal DoF.

\vspace{-0.75em}\subsection{Learnable Bandwidth and Matrix Caching}
\label{app:learnable_results}

\begin{figure}[ht]
\centering
\includegraphics[width=\linewidth]{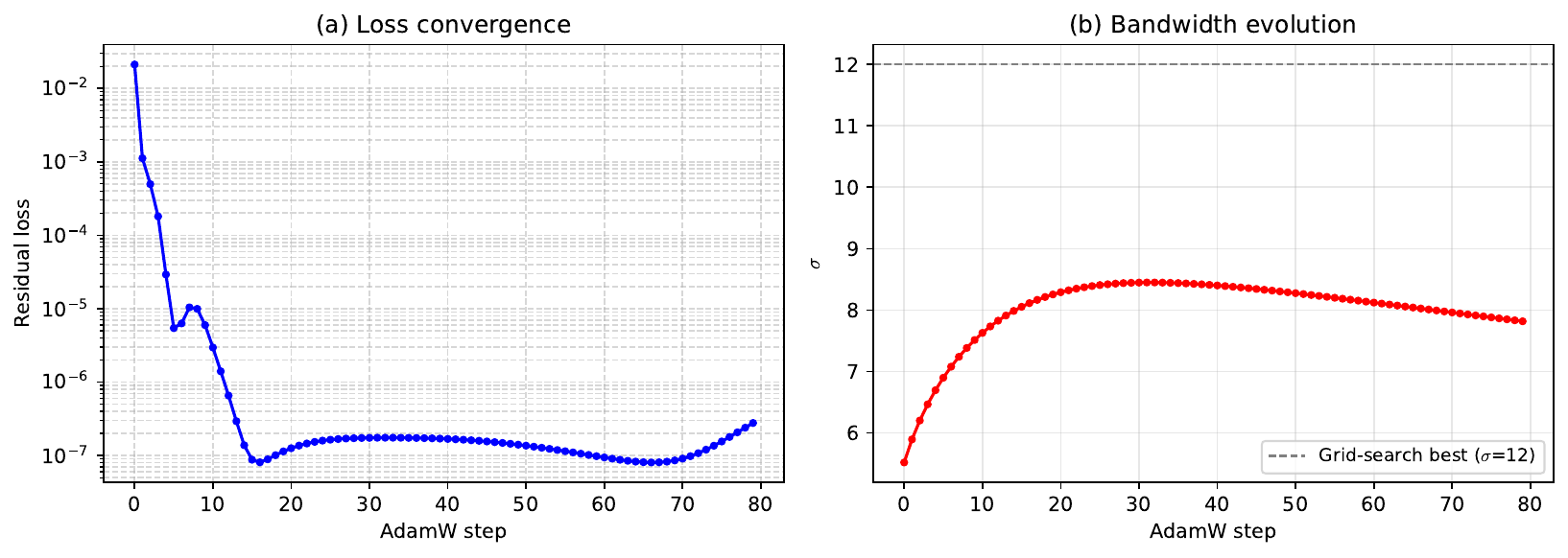}
\caption{Learnable bandwidth: (a) training loss; (b) $\sigma$ evolution (dashed = grid-search optimum).}
\label{fig:learnable_bw}
\vskip -0.1in
\end{figure}

\vspace{-0.75em}\paragraph{Learnable bandwidth.}
On Helmholtz 2D ($k{=}10$), a \texttt{LearnableFastLSQ} module ($N{=}500$, $\sigma_0{=}5$) is optimised via 80 AdamW steps ($\eta{=}0.1$).
At each step, the PDE matrix is assembled analytically and $\bbeta^*(\sigma) = \bA(\sigma)^\dagger\bb$ is solved exactly; only $\sigma$ receives a gradient.
The training loss drops from $2.1\times10^{-2}$ to $8.1\times10^{-8}$ while $\sigma$ evolves from 5 toward the grid-search optimum at 12 (Figure~\ref{fig:learnable_bw}).

\vspace{-0.75em}\paragraph{Matrix caching.}
\label{app:caching_results}
Pre-computing $\bA^\dagger$ once on Helmholtz 2D ($N{=}1500$, $M{=}6000$) reduces each subsequent solve from 453\,ms (full \texttt{lstsq}) to 1.3\,ms (cached multiply), a $362\times$ speedup enabling sub-millisecond parametric sweeps.

\end{document}